# Low Lying Zeros of Families of $L$-Functions

H.Iwaniec, W.Luo and P.Sarnak [1]

## 1 Introduction

In Iwaniec-Sarnak [IS] the percentages of nonvanishing of central values of families of $GL_2$ automorphic L-functions was investigated. In this paper we examine the distribution of zeros which are at or near $s = \frac{1}{2}$ (that is the central point) for such families of L-functions. Unlike [IS], most of the results in this paper are conditional, depending on the Generalized Riemann Hypothesis (GRH). It is by no means obvious, but on the other hand not surprising, that this allows us to obtain sharper results on nonvanishing.

The density and the distribution of zeros near $s = \frac{1}{2}$ for certain families $\mathcal{F}$ of L-functions has been studied recently in Katz-Sarnak [KS1, KS2]. The philosophy and conjectures which emerge assert that for such families of L-functions, the distributions of the low-lying zeros, when we order the L-functions by their conductors (see below), are governed by a symmetry group $G(\mathcal{F})$ associated with $\mathcal{F}$. In the cases where one can identify the function field analogues and compute the scaling limits of the corresponding monodromies of the family, one arrives at such a symmetry $G(\mathcal{F})$. Examples where this can be done and where the corresponding predictions can be verified are given in [KS2]. One of our aims in this paper is to pursue these conjectures

[1] Research partially supported by NSF grants DMS-94-01571, DMS-98-01642, and by AIM.

for families of $GL_2$ and $GL_3$ L-functions.

The families which we consider are the following:

(1). $\mathcal{F}_K$: $L(s,\ f)$'s, where $f$ is a holomorphic cusp form of even integral weight $k \leq K$ for $\mathbf{\Gamma} = SL_2(\mathbf{Z})$, as $K \to \infty$. We define the conductor $c_f$ of $L(s,\ f)$ to be $k^2$.

(2). $\mathcal{F}_N$: $L(s,\ f)$'s, where $f$ is a holomorphic newform of fixed weight $k$ for $\mathbf{\Gamma}_0(N)$ (the Hecke congruence subgroup). Here we let $N \to \infty$, and for simplicity we assume further that $N$ is prime. The conductor $c_f$ of $L(s,\ f)$ in this case is equal to $N$.

(3). $\mathcal{F}_{\mathrm{sym}^2}$: $L(s,\ \mathrm{sym}^2(f))$'s, where $f$ is as in (1), and $\mathrm{sym}^2(f)$ is the symmetric square representation associated to $f$, (see [Shi] for the analytic properties of these L-functions). Define the conductor $c_{\mathrm{sym}^2(f)}$ of $\mathrm{sym}^2(f)$ to be $k^2$.

The families $\mathcal{F}_K$ and $\mathcal{F}_N$ above consist of $GL_2$ forms of general type. Considerations similar to those in [KS2] suggest that their symmetries $G(\mathcal{F}_K)$ and $G(\mathcal{F}_N)$ are $O(\infty)$ (i.e. the scaling limits of orthogonal groups). In particular for each such family approximately half of the self-dual L-functions $L(s,\ f)$ have even functional equations (i.e. with sign $\epsilon_f = 1$) and the other half have odd functional equations ($\epsilon_f = -1$). It is natural to consider the subfamilies $\mathcal{F}_K^{\pm}$ and $\mathcal{F}_N^{\pm}$, where $+$ denotes those $f$'s with $\epsilon_f = 1$ and $-$ those with $\epsilon_f = -1$. The corresponding symmetries are $SO(\text{even})$ (or $O^-(\text{odd})$) and $SO(\text{odd})$ (or $O^-(\text{even})$) respectively.

The family $\mathcal{F}_{\mathrm{sym}^2(f)}$ consists of self-dual $GL_3$ forms (after an application of Gelbart-Jacquet lifting [GJ]). The expected symmetry group for these is $G(\mathcal{F}_{\mathrm{sym}^2}) = Sp(\infty)$. Note that the sign of the functional equation, $\epsilon_{\mathrm{sym}^2(f)}$, is always $+1$.

We investigate the '1-level' densities (see [KS1]) of the low lying zeros. For $F$ any of the $f$'s above or $\text{sym}^2(f)$, denote the nontrivial zeros of $L(s, F)$ by

$$\rho_F = \frac{1}{2} + i\gamma_F \tag{1}$$

Thus $\gamma_F \in \mathbf{R}$ is equivalent to the GRH for $L(s, F)$ being valid. Throughout the paper, let $\phi \in \mathcal{S}(\mathbf{R})$ (a Schwartz function) be even for which the Fourier transform

$$\hat{\phi}(\xi) = \int_{-\infty}^{\infty} \phi(x) e^{-2\pi i x \xi} dx$$

is of compact support (so that $\phi$ extends to an entire function). We define the density sum for zeros near $s = \frac{1}{2}$ by

$$D(F, \phi) = \sum_{\gamma_F} \phi(\gamma_F (\log c_F)/2\pi). \tag{2}$$

Since $\phi$ is localized, the scaling by $(\log c_F)/2\pi$ means that $D$ measures the density of zeros of $L(s, F)$ which are within $1/\log c_F$ of $\frac{1}{2}$. For a given family $\mathcal{F}$, let $M_x(\mathcal{F}) = \#\{F \in \mathcal{F}, c_F \leq x\}$. The density over $\mathcal{F}$ of low lying zeros (ordered by conductor) is defined to be

$$\frac{1}{M_x(\mathcal{F})} \sum_{F \in \mathcal{F}, c_F \leq x} D(F, \phi). \tag{3}$$

This should converge to one of the densities corresponding to each symmetry type as computed in [KS1]. These densities are given by

$$\begin{aligned}
W(SO(\text{even}))(x)dx &= \left(1 + \frac{\sin 2\pi x}{2\pi x}\right) dx, \\
W(SO(\text{odd}))(x)dx &= \delta_0(x) + \left(1 - \frac{\sin 2\pi x}{2\pi x}\right) dx, \\
W(O)(x)dx &= \frac{1}{2}\delta_0 + dx, \\
W(Sp)(x)dx &= \left(1 - \frac{\sin 2\pi x}{2\pi x}\right) dx.
\end{aligned} \tag{4}$$

Here $\delta_0$ is the unit point measure at 0.

The 'Density Conjecture' for any of the families above is the statement

$$\lim_{x \to \infty} \frac{1}{M_x(\mathcal{F})} \sum_{F \in \mathcal{F}, c_F \leq x} D(F, \phi) = \int_{-\infty}^{\infty} \phi(x) W(G(\mathcal{F}))(x) dx \tag{5}$$

for any $\phi \in \mathcal{S}(\mathbf{R})$ with support of $\hat{\phi}$ compact, where $G(\mathcal{F})$ is the proposed symmetry.

Towards this Density Conjecture the following theorem is our main result. Let $H_k^*(\mathbf{\Gamma}_0(N))$ denote the set of newforms (normalized to have first coefficient 1) of weight k for $\mathbf{\Gamma}_0(N)$. Let $H_k^{\pm}(\mathbf{\Gamma}_0(N))$ be the subsets with $\epsilon_f = 1$ and $\epsilon_f = -1$ respectively. Note that if $N = 1$ for $f \in H_k^*(\Gamma)$ we have $\epsilon_f = (-1)^{k/2}$. Set

$$\begin{aligned}
M^+(K) &= \sum_{k \leq K, k \equiv 0(\text{ mod } 4)} \#H_k^*(\mathbf{\Gamma}), \\
M^-(K) &= \sum_{k \leq K, k \equiv 2(\text{ mod } 4)} \#H_k^*(\mathbf{\Gamma}), \\
M(K) &= M^+(K) + M^-(K), \\
M_k^+(N) &= \#H_k^+(\mathbf{\Gamma}_0(N)), \\
M_k^-(N) &= \#H_k^-(\mathbf{\Gamma}_0(N)), \\
M_k^*(N) &= M_k^+(N) + M_k^-(N) = \#H_k^*(\mathbf{\Gamma}_0(N)).
\end{aligned}$$

**Theorem 1**: Assume the GRH for all relevant L-functions. For any $\phi \in \mathcal{S}(\mathbf{R})$ with the support of $\hat{\phi}$ in $(-2, 2)$ we have:

$$\lim_{K \to \infty} \frac{1}{M^+(K)} \sum_{k \equiv 0(\text{ mod } 4), k \leq K} \sum_{f \in H_k^+(\mathbf{\Gamma})} D(f, \phi) = \int_{-\infty}^{\infty} \phi(x) W(SO(\text{even}))(x) dx,$$

$$\lim_{K \to \infty} \frac{1}{M^-(K)} \sum_{k \equiv 2(\text{ mod } 4), k \leq K} \sum_{f \in H_k^-(\mathbf{\Gamma})} D(f, \phi) = \int_{-\infty}^{\infty} \phi(x) W(SO(\text{odd}))(x) dx,$$

$$\lim_{N \to \infty} \frac{1}{M_k^+(N)} \sum_{f \in H_k^+(N)} D(f, \phi) = \int_{-\infty}^{\infty} \phi(x) W(SO(\text{even}))(x) dx,$$

$$\lim_{N \to \infty} \frac{1}{M_k^-(N)} \sum_{f \in H_k^-(N)} D(f, \phi) = \int_{-\infty}^{\infty} \phi(x) W(SO(\text{odd}))(x) dx.$$

Moreover, for any $\phi \in \mathcal{S}(\mathbf{R})$ with the support of $\hat{\phi}$ in $(-4/3,\ 4/3)$, we have:

$$\lim_{K \to \infty} \frac{1}{M(K)} \sum_{k \equiv 0 (\text{ mod } 2), k \leq K} \sum_{f \in H_k(\Gamma)} D(\text{sym}^2(f),\ \phi) = \int_{-\infty}^{\infty} \phi(x) W(Sp)(x) dx.$$

**Remarks**:

(A). Theorem 1 confirms the density conjecture (5) for the families $\mathcal{F}_K, \mathcal{F}_N$ and $\mathcal{F}_{\text{sym}^2}$ in the ranges of the support of $\hat{\phi}$, in which it applies.

(B). Theorem 1 can be established unconditionally for $\phi$'s with the support of $\hat{\phi}$ in $(-1, 1)$. The extension to $\hat{\phi}$'s with larger support is significant. The reason is that the functions $\hat{W}(Sp)(\xi), \hat{W}(SO(\text{even}))(\xi)$ and $\hat{W}(SO(\text{odd}))(\xi)$ all have discontinuities at $\xi = \pm 1$. These discontinuities signal that there are fundamentally new contributions to the asymptotics for such $\phi$'s. Indeed for support of $\hat{\phi}$ in $(-1, 1)$ the contribution comes from the 'diagonal' term in Petersson formula [Pet]. The new non-diagonal terms contributing to the asymptotics arise from Kloosterman sums (such contributions occur in [IS] and earlier in [DFI]). These terms are of a different nature and their appearance is quite subtle for the family $\mathcal{F}_{\text{sym}^2}$. We note that the results in Theorem 1 go well beyond the similar analysis of the pair and higher correlations for the zeros of the Riemann zeta function (Montgomery [Mon], Hejhal [He], Rudnick-Sarnak [RS]). The analysis in those works extends only as far as the diagonal terms being the main contribution to the asymptotics. In particular in as much as Theorem 1 tests the density conjecture (5) beyond the diagonal, we feel it lends strong evidence to the truth of the conjecture.

By choosing positive test functions $\phi$ in Theorem 1 with appropriate restrictions on the support of $\hat{\phi}$, we obtain estimates for the number of $f \in \mathcal{F}$ such that $L(\frac{1}{2}, f) \neq 0$. An analysis of the optimal choice for $\phi$ for this purpose involves extremizing a quadratic form subject to a linear constraint and is carried out in the Appendix.

**Corollary 2**: Assume the GRH for $GL_2$ automorphic L-functions as well as Dirichlet L-functions. Then for the family $\mathcal{F}_K$ we have:

(a) $\liminf\limits_{K\to\infty} \dfrac{\#\{f \in H_k(\Gamma); k \leq K, k \equiv 0(\mathrm{mod}\,4), L(\frac{1}{2}, f) \neq 0\}}{\#\{f \in H_k(\Gamma); k \leq K, k \equiv 0(\mathrm{mod}\,4)\}} > 9/16,$

(b) $\liminf\limits_{K\to\infty} \dfrac{\#\{f \in H_k(\Gamma); k \leq K, k \equiv 2(\mathrm{mod}\,4), L'(\frac{1}{2}, f) \neq 0\}}{\#\{f \in H_k(\Gamma); k \leq K, k \equiv 2(\mathrm{mod}\,4)\}} > 15/16.$

For the family $\mathcal{F}_N$:

(c) $\liminf\limits_{N\to\infty} \dfrac{\#\{f \in H_k^+(\Gamma_0(N)); L(\frac{1}{2}, f) \neq 0\}}{\#\{f \in H_k^+(\Gamma_0(N))\}} > 9/16,$

(d) $\liminf\limits_{N\to\infty} \dfrac{\#\{f \in H_k^-(\Gamma_0(N)); L'(\frac{1}{2}, f) \neq 0\}}{\#\{f \in H_k^-(\Gamma_0(N))\}} > 15/16,$

(e) $\limsup\limits_{N\to\infty} \dfrac{\sum_{f \in H_k^*(N)} \mathrm{ord}(\frac{1}{2}, f)}{M_k^*(N)} < 1,$

where $\mathrm{ord}(s_0, f)$ is the order of vanishing at $s = s_0$ of $L(s, f)$.

For the family $\mathcal{F}_{\mathrm{sym}^2}$:

(f) $\liminf\limits_{K\to\infty} \dfrac{\#\{f \in H_k(\Gamma); k \leq K, k \equiv 0(\mathrm{mod}\,2), L(\frac{1}{2}, \mathrm{sym}^2(f)) \neq 0\}}{\#\{f \in H_k(\Gamma); k \leq K, k \equiv 0(\mathrm{mod}\,2)\}} > 55/64.$

**Remarks**:

(A). The Density Conjecture would yield the (presumably) true values for the above limits, that is the value of 1 for all but case (e) and the value $\frac{1}{2}$ for (e).

(B). All the bounds in Corollary 2 depend crucially on the extensions of the range of the support of $\hat{\phi}$ in Theorem 1. Without these extensions one would for example obtain a lower bound of greater or equal to 1/2 in place of 9/16 in (a) and (c). In [IS] this lower bound of greater or equal to 1/2 is established unconditionally. Moreover it is shown that improving this to anything bigger than 1/2 is intimately connected to the Landau-Siegel zero. Of course, in Corollary 2 this is not an issue since we are assuming GRH. The minute improvements implied by the strict inequalities in the Corollary can be determined explicitly (see Appendix 1).

(C). An unconditional lower bound of greater than or equal to 7/8 for (d) above has recently been established by Kowalski-Michel [KM] and Vanderkam [Va]. Also in Kowalski-Michel [KM] an unconditional upper bound of 6.5 in place of 1 in (e) above, is proven.

(D). The Density Conjecture for $\mathcal{F}_N$ and the cases (c), (d) and (e) of Corollary 2 specified to $k = 2$ have applications to the estimation of the ranks of the Jacobians of the modular curves $X_0(N)$, see [KS2] for a description of these implications.

The results in Theorem 1 are stated for families with fixed sign of $\epsilon_f$, and they are conditional. However, if we do not break parity (with respect to the sign of $\epsilon_f$) we can easily establish unconditionally a version of Theorem 1, which is as strong . We do so for a weighted version of the densities for the family $\mathcal{F}_\mathcal{K}$, the weight being $L(1, \text{sym}^2(f))^{-1}$ which is dictated by the Petersson formula [Pet].

Fix $h \in C_0^\infty(0, \infty)$, $h \geqslant 0$ and define the averaging operator $A_K$ by

$$A_K[X_f] = \sum_{k \equiv 0(\bmod 2)} h(k/K) \sum_{f \in H_k(\mathbf{\Gamma})} L(1, \text{sym}^2(f))^{-1} X_f, \tag{6}$$

where $X$ is a function from $H_k(\mathbf{\Gamma})$ to $\mathbf{R}$.

**Theorem 3**: For $\phi \in \mathcal{S}(\mathbf{R})$ with support$(\hat{\phi}) \subset (-2, 2)$, and $D_k(f, \phi)$ the approximation to $D(f, \phi)$ given in (20), we have

$$\lim_{K \to \infty} \frac{A_K[D(f, \phi)]}{A_K[1]} = \int_{-\infty}^\infty \phi(x) W(O)(x) dx.$$

Under some natural hypothesis about exponential sums of analytic (as opposed to arithmetic) functions over primes in progressions we can extend the range $(-2, 2)$ and this has striking consequences.

**Hypothesis 4**: For any $\epsilon > 0$, $X \geq 2$, $c \geq 1$ and $a$ with $(a, c) = 1$ we have

$$\sum_{p \leq X, p \equiv a (\bmod c)} e(2\sqrt{p}/c) \ll_\epsilon X^{\frac{1}{2}+\epsilon}.$$

where the implied constant depends only on $\epsilon$.

A discussion of this Hypothesis as well as its relation to the zeros of Dirichlet L-functions (i.e. $GL_1/\mathbf{Q}$ L-functions) is given in [CIRS]. Here we note that for our purpose of extending the range $(-2, 2)$, all that we need from Hypothesis 4 is a bound of the form $X^\delta c^A$ with some $\delta < 3/4$ and some $A > 0$. With $\delta = 7/8$ such a bound is essentially due to Vinogradov [V] while the bound with $\delta = 3/4$ follows from the Density Hypothesis for Dirichlet L-functions.

**Theorem 5**: Assume Hypothesis 4, then Theorem 3 is valid for $\phi$'s for which $\hat{\phi}(\xi)$ is supported in $(-7/3, 7/3)$.

The remarkable thing about Theorem 5 (or any improvement over the range $(-2, 2)$) is that it strikes at the $GL_2$ Riemann Hypothsis! That is, it implies

**Corollary 6**: Assume Hypothesis 4 and that the zeros of $L(s, f)$ for $f \in H_k^*(\Gamma)$ are either on $\Re(s) = \frac{1}{2}$ or on the real axis. Then for all $k$ sufficiently large (and effective) and any $f \in H_k^*(\Gamma)$ we have

$$L(\beta, f) \neq 0, \text{if } \beta > 13/14.$$

**Remark**:

By a somewhat different method we can remove the annoying assumption about $L(s, f)$ having its zeros either on $\Re(s) = \frac{1}{2}$ or being real [ILS]. This leads to Hypothesis 4 implying a quasi Riemann Hypothesis for $L(s, f)$ for $f$ of large weight $k$. What is remarkable is that Hypothesis 4 is concerned with bounds for exponential sums of an analytic function over primes in progressions (so it is a classical $GL_1$ type assumption) which implies a quasi Riemann Hypothesis for the $GL_2$ L-functions(it reminds one of the Lang-Weil Theorem [LW] which gives a quasi Riemann Hypothesis for zeta function of varieties over finite fields by using RH for zeta functions of curves). It is very different to the familiar relation between cancellations in $\sum_{p \leq X} \lambda_f(p)$, where $\lambda_f(p)$ is the $p$-th Fourier coefficient of $f$, and the Riemann Hypothesis for $L(s, f)$. The relation in Corollary 6 also throws new light on these classical exponential sums in Hypothesis 4.

## 2 Large weight

### 2.1 Preliminaries

In order to carry out summations over the weight $k$ ($k$ even) we recall some results about series of Bessel functions, the Petersson Formula and the Explicit Formula.

**Proposition 1**. Fix a real valued function $h \in C_0^\infty(0, \infty)$. For $a = \pm 1$ and $x$, $L > 0$ we have

(A).

$$4 \sum_{l \equiv a (\bmod 4)} h(\frac{l}{L}) J_l(x) = h\left(\frac{x}{L}\right) + i^{1-a} \Im \left( e^{i(x-\pi/4)} \frac{L}{\sqrt{x}} V_h\left(\frac{L^2}{2x}\right) \right) + \frac{x}{6L^3} h^{(3)}\left(\frac{x}{L}\right) + O\left(\frac{x}{L^4} + \frac{x^2}{L^6}\right).$$

(B).

$$2 \sum_{l \equiv 1 (\bmod 2)} h(\frac{l}{L}) i^{l+1} J_l(x) = -\Im \left( e^{i(x-\pi/4)} \frac{L}{\sqrt{x}} V_h\left(\frac{L^2}{2x}\right) \right) + O\left(\frac{x}{L^4} + \frac{x^2}{L^6}\right),$$

where $V_h \in \mathcal{S}(\mathbf{R})$ and is given by

$$V_h(\xi) = \frac{1}{\sqrt{2\pi}} \int_0^\infty \frac{h(\sqrt{u})}{\sqrt{u}} e^{i\xi u} du.$$

**Proof.** This follows easily from the integral representation

$$J_l(x) := \int_{-\frac{1}{2}}^{\frac{1}{2}} e(lt)\exp(-x\sin(2\pi t))dt,$$

and Poisson summation, see [Iw1, pp 85], where for $g(x)$ we are using the function $h(\frac{x}{L})$. □

Note, in particular, that since $V_h$ is rapidly decreasing at $\infty$ the imaginary part term in Proposition 1 is $O(L^{-A})$ for any $A > 0$, if $x \leq L^{2-\delta}$ for some $\delta > 0$.

Recall that $H_k^*(\Gamma), k$ even, denotes the set of holomorphic Hecke cusp forms (with first Fourier coefficient $a_1 = 1$) of weight $k$ for $\Gamma = SL_2(\mathbf{Z})$. As is well known

$$\#H_k(\Gamma) = \begin{cases} [k/12] - 1 & \text{if } k \equiv 2 (\text{mod } 12); \\ [k/12] & \text{if otherwise.} \end{cases} \qquad (7)$$

For such an $f$ let $L(s, f)$ be its $L$-function:

$$L(s, f) = \sum_{n=1}^{\infty} \frac{\lambda_f(n)}{n^s}, \qquad (8)$$

where $\lambda_f(n)$ are the normalized Hecke eigenvalues of $f$. In particular, $\lambda_f(n) \in \mathbf{R}$ and it satisfies the following multiplication rules

$$\lambda_f(n)\lambda_f(m) = \sum_{d|(n,m)} \lambda_f\left(\frac{mn}{d^2}\right). \qquad (9)$$

Moreover we have ([De])

$$|\lambda_f(n)| \leq \tau(n) := \sum_{d|n} 1. \qquad (10)$$

The functions

$$\Lambda(s, f) = (2\pi)^{-s}\Gamma(s + (k-1)/2)L(s, f) \qquad (11)$$

are entire and satisfy the self-dual functional equations

$$\Lambda(s, f) = i^k \Lambda(1-s, f). \qquad (12)$$

The symmetric square $L$-function $L(s, \text{sym}^2(f))$ corresponding to $f$ is defined by

$$L(s, \text{sym}^2(f)) = \zeta(2s) \sum_{n=1}^{\infty} \frac{\lambda_f(n^2)}{n^s}. \tag{13}$$

Shimura ([Sh]) has shown that this Euler product of degree 3 is entire and satisfies a functional equation

$$\begin{aligned}\Lambda(s, \text{sym}^2(f)) &:= (\pi)^{-3s/2} \Gamma\left(\frac{s+1}{2}\right) \Gamma\left(\frac{s+k-1}{2}\right) \Gamma\left(\frac{s+k}{2}\right) L(s, \text{sym}^2(f)) \\ &= \Lambda(1-s, \text{sym}^2(f)).\end{aligned} \tag{14}$$

It is known that $L(s, \text{sym}^2(f)) \neq 0$ for $\Re(s) = 1$ and also ([HL]) that

$$(\log k)^{-2} \ll L(1, \text{sym}^2(f)) \ll (\log k)^2. \tag{15}$$

The key tool for averaging over forms of weight $k$ is the following

**Proposition 2** ([Pet]).

$$\frac{2\pi^2}{k-1} \sum_{f \in H_k^*(\Gamma)} \frac{\lambda_f(m)\lambda_f(n)}{L(1, \text{sym}^2(f))} = \delta(m, n) + 2\pi i^k \sum_{c=1}^{\infty} \frac{S(m, n; c)}{c} J_{k-1}\left(\frac{4\pi\sqrt{mn}}{c}\right),$$

where $\delta(m, n)$ is 1 if $m = n$ and 0 otherwise, while $S(m, n; c)$ is the Kloosterman sum

$$S(m, n; c) = \sum_{x(\bmod c), x\overline{x} \equiv 1(\bmod c)} e\left(\frac{mx + n\overline{x}}{c}\right).$$

Next we define averaging operators $A_K$ and $A_K^{\pm}$ as follows: For $h \geq 0$ in $C_0^{\infty}(0, \infty)$ and $X_f$ a quantity depending on $f$, set

$$A_K[X_f] = \sum_{k \equiv 0(\bmod 2)} h\left(\frac{k-1}{K}\right) \frac{2\pi^2}{k-1} \sum_{f \in H_k^*(\Gamma)} \frac{X_f}{L(1, \text{sym}^2(f))} \tag{16}$$

and

$$A_K^{\epsilon}[X_f] = \sum_{k \equiv 1-\epsilon(\bmod 4)} h\left(\frac{k-1}{K}\right) \frac{2\pi^2}{k-1} \sum_{f \in H_k^*(\Gamma)} \frac{X_f}{L(1, \text{sym}^2(f))}, \tag{17}$$

where $\epsilon = \pm 1$. Thus according to (12) $A_K^\epsilon$ averages over $f$'s whose sign of their functional equation is $\epsilon$.

Applying Propositions 1 and 2 with $A_K[1] = A_K[\lambda_f^2(1)]$ we get

$$A_K[1] = \frac{1}{2}\hat{h}(0)K + O(K^{-4}), \tag{18}$$

and

$$A_K^{\pm 1}[1] = \frac{1}{4}\hat{h}(0)K + O(K^{-4}). \tag{19}$$

In particular the average of $L(1, \text{sym}^2(f))^{-1}$ is $6/\pi^2$ which is best thought of as $1/\zeta(2)$.

The final result which we will use repeatedly is the explicit formula relating the zeros of $\Lambda(s, f)$, which we denote by $\rho_f = \frac{1}{2} + i\gamma_f$, and the numbers $\lambda_f(n)$. Write $\lambda_f(p) = \alpha_f(p) + \beta_f(p)$ with $\alpha_f(p)\beta_f(p) = 1$ (so that according to (10), $|\alpha_f(p)| = |\beta_f(p)| = 1$).

**Proposition 3** ([RS]). Let $\phi \in \mathcal{S}(\mathbf{R})$ be even and for which $\hat{\phi}(\xi) = \int_{-\infty}^{\infty} e(-x\xi)\phi(x)dx$ is of compact support. Then for $f \in H_k^*(\Gamma)$

$$\sum_{\gamma_f} \phi\left(\frac{\gamma_f}{2\pi}\right) = \frac{1}{2\pi}\int_{-\infty}^{\infty}\left(-\log(2\pi) + \frac{\Gamma'}{\Gamma}\left(\frac{k+1}{2} + ir\right) + \frac{\Gamma'}{\Gamma}\left(\frac{k+1}{2} - ir\right)\right)\phi\left(\frac{r}{2\pi}\right)dr$$
$$- 2\sum_p \sum_{\nu=1}^{\infty} \frac{\log p}{p^{\nu/2}}(\alpha_f(p)^\nu + \beta_f(p)^\nu)\hat{\phi}(\log p^\nu).$$

## 2.2  When $k \equiv 0 (\text{mod} 2)$:

Fix a test function $\phi$ as in Proposition 3. For $K$ large and $A_1 K \leq k \leq A_2 K$ ($A_1 < A_2$, fixed) define the densities of low-lying zeros for $f \in H_k^*(\Gamma)$ by

$$D_K(f, \phi) = \sum_{\gamma_f} \phi\left(\frac{\gamma_f \log K^2}{2\pi}\right). \tag{20}$$

These serve as approximations to the densities (see (2))

$$D(f, \phi) := \sum_{\gamma_f} \phi\left(\frac{\gamma_f \log k^2}{2\pi}\right). \tag{21}$$

We examine the asymptotic behavior of averages $A_K[D_K(f, \phi)]$ as $K \to \infty$. Applying the explicit formula (Proposition 3) to the test function $\phi((2\log K)x)$ whose Fourier transform is $\frac{1}{2\log K}\hat{\phi}\left(\frac{\xi}{2\log K}\right)$, together with the relation $\alpha_f^2(p) + \beta_f^2(p) = \lambda_f(p^2) - 1$, one gets

$$\begin{aligned}
\sum_{\gamma_f} \phi\left(\frac{\gamma_f \log K^2}{2\pi}\right) &= \int_{-\infty}^{\infty} \phi(x)dx + \frac{1}{\log K}\sum_p \frac{\log p}{p}\hat{\phi}\left(\frac{\log p}{\log K}\right) \\
&\quad - \frac{1}{\log K}\sum_p \frac{\log p}{\sqrt{p}}\lambda_f(p)\hat{\phi}\left(\frac{\log p}{2\log K}\right) \\
&\quad - \frac{1}{\log K}\sum_p \frac{\log p}{p}\lambda_f(p^2)\hat{\phi}\left(\frac{\log p}{\log K}\right) + O\left(\frac{1}{\log K}\right). \tag{22}
\end{aligned}$$

Hence $A_K[D_K(f, \phi)]$ splits accordingly into

$$A_K[D_K(f, \phi)] = I + II + III + IV + O\left(\frac{A_K[1]}{\log K}\right). \tag{23}$$

Clearly

$$I = A_K\left[\int_{-\infty}^{\infty} \phi(x)dx\right] = \left(\int_{-\infty}^{\infty} \phi(x)dx\right)A_K[1], \tag{24}$$

while

$$II = A_K[1]\frac{1}{\log K}\sum_p \frac{\log p}{p}\hat{\phi}\left(\frac{\log p}{\log K}\right)$$

which by the prime number theorem is asymptotically equal to

$$A_K[1]\frac{1}{\log K}\int_0^{\infty} \hat{\phi}\left(\frac{\log t}{\log K}\right)\frac{dt}{t} = A_K[1]\int_0^{\infty} \hat{\phi}(y)dy = \frac{\phi(0)}{2}A_K[1]. \tag{25}$$

The assumption that support $(\hat{\phi}) \subset (-2, 2)$ restricts the ranges of summation in III to $p \leq K^{4-\delta}$ for some $\delta > 0$ and in IV to $p \leq K^{2-\delta}$. Thus applying Proposition 2 and part B of Proposition 1 (with $L = K$, $x = 4\pi c^{-1}\sqrt{p} \ll c^{-1}K^{2-\delta/2}$, so that $L^2/2x \gg cK^{\delta/2}$ and the first term on the right-hand side of (B) is absorbed by the error term $xL^{-4}$) we get

$$III \ll \frac{1}{\log K} \sum_{p \leq K^{4-\delta}} \frac{\log p}{\sqrt{p}} \sum_c \frac{|S(1,p;c)|}{c} \frac{\sqrt{p}}{cK^4}.$$

It is convenient here and elsewhere to invoke the well-known bound of Weil [Wei]

$$|S(m,n;c)| \leq (m,n,c)^{\frac{1}{2}} \tau(c) c^{\frac{1}{2}}. \qquad (26)$$

This yields

$$III \ll K^{-\delta}. \qquad (27)$$

Similarly we may show that

$$IV \ll \frac{1}{\log K} \sum_{p \leq K^{2-\delta}} \frac{\log p}{p} \sum_c \frac{|S(1,p^2;c)|}{c} \frac{p}{cK^4} \ll K^{-2}. \qquad (28)$$

Combining these estimates we conclude that if support $(\hat{\phi}) \subset (-2, 2)$, then

$$A_K[D_K(f,\phi)] \sim (\int_{-\infty}^{\infty} \phi(x)dx + \frac{\phi(0)}{2})A_K[1]$$

$$\sim (\int_{-\infty}^{\infty} \phi(x)dx + \frac{\phi(0)}{2})\frac{\hat{h}(0)K}{2}, \qquad (29)$$

as $K \to \infty$. This establishes Theorem 3.

We turn to extending this result to $\hat{\phi}$'s with larger support. Allowing $p$ to be larger than $K^4$ in the term $III$ above one faces the new term coming from Proposition 1 part B, viz the imaginary part of

$$\sum_p \frac{\log p}{\sqrt{p}} \hat{\phi}\left(\frac{\log p}{2\log K}\right) \sum_c \frac{S(1,p,c)}{c}$$

$$\times \left(e^{i(4\pi\sqrt{p}/c - \pi/4)} \frac{K\sqrt{c}}{\sqrt{4\pi}p^{1/4}} V_h\left(\frac{K^2 c}{8\pi\sqrt{p}}\right) + O\left(\frac{\sqrt{p}}{cK^4}\right) + O\left(\frac{p}{c^2 K^6}\right)\right). \qquad (30)$$

Suppose that the support of $\hat{\phi}$ imposes the restriction $p \leq P$, then the O-term in (30) contributes (here and henceforth, $\epsilon$ is an arbitrarily small positive number and it may change its value at each occurence)

$$\ll PK^{-4+\epsilon} + P^{3/2}K^{-6+\epsilon} = o(K), \tag{31}$$

provided that $P < K^{14/3-\delta}$, i.e. if

$$\text{support } \hat{\phi} \subset (-7/3,\ 7/3). \tag{32}$$

As to the other term in (30), it may be written as

$$\sum_{c \leq P^{\frac{1}{2}+\epsilon}/K^2} \frac{1}{\sqrt{c}} \sum_{a(\bmod c)}^{*} S(1,a;c) \int_0^\infty \hat{\phi}\left(\frac{\log t}{2 \log K}\right) V_h\left(\frac{cK^2}{8\pi\sqrt{t}}\right) \cdot t^{-3/4} dF_{a,c}(t) + O(1), \tag{33}$$

where

$$F_{a,c}(t) = \sum_{p \equiv a(\bmod c), p \leq t} (\log p) e(2\sqrt{p}/c). \tag{34}$$

The restriction on the $c$ sum follows from $V_h$ being in $\mathcal{S}(\mathbf{R})$, and of course if $P \leq K^4$ then the term (33) is essentially not present. According to Hypothesis 4 of the Introduction which we now assume we have

$$F_{a,c}(t) \ll_\epsilon t^{\frac{1}{2}+\epsilon}. \tag{35}$$

Integrating by parts in (33) and invoking (35) yields the bound of $K^{-3}P^{\frac{3}{4}+\epsilon}$ for the quantity in (33). Again this is $o(K)$ as long as support $\hat{\phi} \subset (-7/3,\ 7/3)$.

Concerning the term IV in (23), an application of Proposition 2 followed by Proposition 1 part B and a direct estimation using (26) shows that

$$IV = o(K)$$

as long as support $\hat{\phi} \subset (-4, 4)$. Combining these estimates we conclude the following Proposition (which is Theorem 5 of the Introduction):

**Proposition 4**. Assume Hypothesis 4 of the Introduction. For $\phi$ with support $\hat{\phi}$ contained in $(-7/3, 7/3)$ (29) holds.

Assume now that $L(s, f)$ satisfies the condition that all its zeros are on $\Re(s) = \frac{1}{2}$ or are real. Then if we choose a test function $\phi(x)$ in the form

$$\phi(x) = \left(\frac{\sin(2\pi\nu x)}{2\pi\nu x}\right)^2, \tag{36}$$

(or even a $C^\infty$ smoothed version of it) where $\nu < 7/6$, then support $\hat{\phi} \subset (-7/3, 7/3)$. Moreover under the assumption on the zeros of $L(s, f)$ we have $D_K(f, \phi) \geq 0$ for any $f$. Applying Proposition 4 we conclude that for any $f$ of weight $k$ (and set $K = k$)

$$\frac{1}{L(1, \text{sym}^2(f))} D(f, \phi) = \sum_{\gamma_f} \left(\frac{\sin(2\nu\gamma_f \log k)}{2\nu\gamma_f \log k}\right)^2 \ll k^2. \tag{37}$$

According to (15) this gives

$$\left(\frac{\sin(2\nu\gamma_f \log k)}{2\nu\gamma_f \log k}\right)^2 \ll_\epsilon k^{2+\epsilon} \tag{38}$$

for any zero $\rho_f = \frac{1}{2} + i\gamma_f$. In particular if $\rho_f \in (\frac{1}{2}, 1)$ so that $\gamma_f = -i\mu_f$, $0 < \mu_f < \frac{1}{2}$, then

$$\left(\frac{\sinh(2\nu\mu_f \log k)}{2\nu\mu_f \log k}\right)^2 \ll_\epsilon k^{2+\epsilon}$$

or

$$k^{4\nu\mu_f} \ll_\epsilon k^{2+\epsilon}. \tag{39}$$

For $k$ large this implies that $\mu_f \leq (2\nu)^{-1}$, i.e. $\mu_f \leq 3/7$. We conclude that any real zero $\rho_f$ satisfies $\rho_f \leq 13/14$. This establishes Corollary 6.

## 2.3 When $k \equiv a \pmod{4}$:

In this subsection we consider similar averages to the previous one except that we break parity in terms of the signs of the functional equations of $L(s, f)$ and we also average with respect to uniform weights. We will also assume GRH - that is the Grand Riemann Hypothesis for all $L$-functions.

Let $h$ and $\phi$ be as before. For $a = 0$ or $2$ consider the sums

$$S_K^{(1-a)}[D_K(f, \phi)] = \sum_{k \equiv a (\bmod 4)} h\left(\frac{k-1}{K}\right) \frac{12}{k-1} \sum_{f \in H_k^*(\Gamma)} D_K(f, \phi) \qquad (40)$$

$$= A_K^{(1-a)}[6\pi^{-2} L(1, \text{sym}^2(f)) D_K(f, \phi)]. \qquad (41)$$

It is easy to see that given $\phi$ as above we can find a $\phi_1 \in \mathcal{S}(\mathbf{R})$ whose Fourier transform has compact support and such that

$$|\phi(Lx)| \leq \phi_1(x), \quad \text{for } L \geq 1, \ x \in \mathbf{R}. \qquad (42)$$

Since we are assuming the Riemann Hypothesis for $L(s, f)$ we have $\gamma_f \in \mathbf{R}$ and hence

$$|D_K(f, \phi)| \leq \sum_{\gamma_f} \phi_1\left(\frac{\gamma_f}{2\pi}\right). \qquad (43)$$

Applying the explicit formula to the latter (and that $k \sim K$) we get

$$D_K(f, \phi) \ll \log K. \qquad (44)$$

The Riemann Hypothesis for $L(s, \text{sym}^2(f))$ allows us to approximate $L(1, \text{sym}^2(f))$ by a short Dirichlet series (in fact one can show unconditionally that such approxima-tions are valid with at most few exceptions (see [IS] and [L])). More precisely shifting the contour integral

$$\int_{\Re(s)=1} L(s+1,\ \mathrm{sym}^2(f))X^s \frac{ds}{s}$$

to $\Re(s) = -\frac{1}{2}$ and applying the Lindelöf Hypothesis (which is a consequence of GRH) for $L(s,\ \mathrm{sym}^2(f))$ yields:

$$L(1,\ \mathrm{sym}^2(f)) = \sum_{m^2 n \leq X} \frac{\lambda_f(n^2)}{m^2 n} + O_\epsilon(k^\epsilon X^{-\frac{1}{2}}). \tag{45}$$

In particular if $\delta_0 > 0$ (very small and fixed) then

$$L(1,\ \mathrm{sym}^2(f)) = \sum_{m^2 n \leq k^{\delta_0}} \frac{\lambda_f(n^2)}{m^2 n} + O(k^{-\delta_0/2}). \tag{46}$$

Moreover from (45) and (15) it follows that for $k$ large

$$\sum_{m^2 n \leq k^{\delta_0}} \frac{\lambda_f(n^2)}{m^2 n} \gg (\log k)^{-2}. \tag{47}$$

Throughout this subsection we will assume that support $\hat{\phi} \subset (-2,\ 2)$ so that for suitable $\delta > 0$

$$\text{support } \hat{\phi} \subset (-2+\delta,\ 2-\delta). \tag{48}$$

Let $\delta_0 = \delta/100$ (say), then using the approximation (46) in (41) and using (44) as well we get

$$S_K^{(1-a)}[D_K(f,\ \phi)] = 6\pi^{-2} A_K^{(1-a)} \left[ D_K(f,\ \phi) \sum_{m^2 n \leq k^{\delta_0}} \frac{\lambda_f(n^2)}{m^2 n} \right] + O(A_K^{(1-a)}[1] K^{-\delta_0/2} \log K). \tag{49}$$

Applying the explicit formula in the form (22) to the middle term in (49) and using (47) we get

$$S_K^{(1-a)}[D_K(f, \phi)]$$

$$= \frac{6}{\pi^2}\left(\int_{-\infty}^{\infty}\phi(x)dx + \frac{1}{\log K}\sum_p \frac{\log p}{p}\hat{\phi}\left(\frac{\log p}{\log K}\right)\right) A_K^{(1-a)}\left[\sum_{m^2 n \leq k^{\delta_0}} \frac{\lambda_f(n^2)}{m^2 n}\right]$$

$$- \frac{6}{\pi^2 \log K}\sum_p \frac{\log p}{\sqrt{p}}\hat{\phi}\left(\frac{\log p}{2\log K}\right) A_K^{(1-a)}\left[\sum_{m^2 n \leq k^{\delta_0}} \frac{\lambda_f(p)\lambda_f(n^2)}{m^2 n}\right]$$

$$- \frac{6}{\pi^2 \log K}\sum_p \frac{\log p}{p}\hat{\phi}\left(\frac{\log p}{\log K}\right) A_K^{(1-a)}\left[\sum_{m^2 n \leq k^{\delta_0}} \frac{\lambda_f(p^2)\lambda_f(n^2)}{m^2 n}\right] +$$

$$+ O\left(\frac{1}{\log K} A_K^{(1-a)}\left[\sum_{m^2 n \leq k^{\delta_0}} \frac{\lambda_f(n^2)}{m^2 n}\right] + K^{1-\delta_0/2}\log K\right). \tag{50}$$

Next we examine the contributions of each of these terms as $K \to \infty$. Applying Propositions 1 and 2 one gets

$$A_K^{(1-a)}\left[\sum_{m^2 n \leq k^{\delta_0}} \frac{\lambda_f(n^2)}{m^2 n}\right] = \frac{K\hat{h}(0)}{4}\sum_{m^2 \leq K^{\delta_0}} \frac{1}{m^2} + O(K^{3\delta_0}K^{3\delta_0 - 4}). \tag{51}$$

Thus the first term in (50) is asymptotic to

$$\frac{K\hat{h}(0)}{4}\left(\int_{-\infty}^{\infty}\phi(x)dx + \frac{\phi(0)}{2}\right) \sim S_K^{(1-a)}[1]\left(\int_{-\infty}^{\infty}\phi(x)dx + \frac{\phi(0)}{2}\right).$$

It also follows from the above that the O-term in (50) is $o(K)$ as $K \to \infty$.

The third term in (50) is also $o(K)$ as $K \to \infty$ since according to Proposition 1 and Proposition 2

$$\sum_p \frac{\log p}{p}\hat{\phi}\left(\frac{\log p}{\log K}\right) A_K^{(1-a)}\left[\sum_{m^2 n \leq k^{\delta_0}} \frac{\lambda_f(p^2)\lambda_f(n^2)}{m^2 n}\right]$$

$$= \sum_{m^2 p \leq K^{\delta_0}} \frac{\log p}{p^2 m^2} \hat{\phi}\left(\frac{\log p}{\log K}\right) \frac{K\hat{h}(0)}{4}$$

$$+ O\left(\sum_p \frac{\log p}{p} \left|\hat{\phi}\left(\frac{\log p}{\log K}\right)\right| \sum_{m^2 n \leq K^{\delta_0}} \frac{1}{m^2 n} \sum_c \frac{|S(p^2, n^2; c)|}{c} \left|h\left(\frac{4\pi pn}{cK}\right)\right|\right)$$

$$\ll K + \sum_{p \leq K^2} \frac{\log p}{p} \sum_{m^2 n \leq K^{\delta_0}} \frac{1}{m^2 n} \sum_{c \leq K} c^{-\frac{1}{2}+\epsilon} \ll K + K^{\frac{1}{2}+\epsilon} \ll K. \tag{52}$$

This leaves us with the second term in (50) which in fact contributes to the main term (this being an 'off-diagonal' contribution). Applying Proposition 2 it equals

$$-\frac{6}{\pi^2 \log K} \sum_p \frac{\log p}{\sqrt{p}} \hat{\phi}\left(\frac{\log p}{2 \log K}\right)$$

$$\times \left[2\pi \sum_{k \equiv a(\bmod 4)} i^k h\left(\frac{k-1}{K}\right) \sum_{m^2 n \leq K^{\delta_0}} \frac{1}{m^2 n} \sum_c \frac{S(p, n^2; c)}{c} J_{k-1}\left(\frac{4\pi n \sqrt{p}}{c}\right)\right].$$

Applying Proposition 1 part A to the above yields (using $p \ll K^{4-2\delta}$)

$$= -\frac{12\pi i^a}{\pi^2 \log K} \sum_p \frac{\log p}{\sqrt{p}} \hat{\phi}\left(\frac{\log p}{2 \log K}\right)$$

$$\times \sum_{m^2 n \leq K^{\delta_0}} \frac{1}{m^2 n} \sum_c \frac{S(p, n^2; c)}{c} \left[\frac{1}{4} h\left(\frac{4\pi n \sqrt{p}}{cK}\right) + O\left(\frac{n\sqrt{p}}{cK^3}\right)\right]$$

$$= -\frac{3i^a}{\pi \log K} \sum_{m^2 n \leq K^{\delta_0}} \frac{1}{m^2 n} \sum_c \frac{1}{c} \sum_p \frac{S(p, n^2; c) \log p}{\sqrt{p}} \hat{\phi}\left(\frac{\log p}{2 \log K}\right) h\left(\frac{4\pi n \sqrt{p}}{cK}\right)$$

$$+ O(K^{-\delta/2})$$

$$= -\frac{3i^a}{\pi \log K} \sum_{m^2 n \leq K^{\delta_0}} \frac{1}{m^2 n} \sum_c \frac{1}{c} \sum_{a(\bmod c)}^{*} S(a, n^2; c)$$

$$\times \sum_{p \equiv a(\bmod c)} \frac{\log p}{\sqrt{p}} \hat{\phi}\left(\frac{\log p}{2 \log K}\right) h\left(\frac{4\pi n \sqrt{p}}{cK}\right) + O(K^{-\delta/2})$$

$$= -\frac{3i^a}{\pi \log K} \sum_{m^2 n \leq K^{\delta_0}} \frac{1}{m^2 n} \sum_c \frac{1}{c} \sum_{a(\bmod c)}^{*} S(a, n^2; c) \frac{1}{\varphi(c)} \sum_{\chi(\bmod c)} \bar{\chi}(a)$$

$$\times \sum_p \frac{\log p}{\sqrt{p}} \chi(p) \hat{\phi}\left(\frac{\log p}{2 \log K}\right) h\left(\frac{4\pi n \sqrt{p}}{cK}\right) + O(K^{-\delta/2}). \tag{53}$$

We next invoke GRH for Dirichlet $L$-functions which gives for $\chi(\mathrm{mod}\, c)$

$$\psi(x,\chi) = \sum_{p \leq x} \chi(p)\log p = \delta_\chi x + O(c^\epsilon x^{\frac{1}{2}+\epsilon}) \tag{54}$$

where $\delta_\chi = 1$ if $\chi$ is principal and is 0 otherwise.

We consider first the contribution to (53) when $\chi \neq \chi_0$ (the principal character). This is

$$\ll \frac{1}{\log K} \sum_{m^2 n \leq K^{\delta_0}} \frac{1}{m^2 n} \sum_{c \leq K^{1-\delta/10}} \frac{1}{c\varphi(c)} \sum_{\chi(\mathrm{mod}\, c), \chi \neq \chi_0} \tag{55}$$

$$\times \left|\sum^*_{a(\mathrm{mod}\, c)} S(a, n^2; c)\bar{\chi}(a)\right| \left|\int_0^\infty \hat{\phi}\left(\frac{\log t}{2\log K}\right) h\left(\frac{4\pi n\sqrt{t}}{cK}\right) \frac{d\psi(t,\chi)}{\sqrt{t}}\right|.$$

Integrating by parts and using (54), together with

$$|\sum^*_{a(\mathrm{mod}\, c)} S(a, n^2; c)\bar{\chi}(a)| \ll_\epsilon c^{1+\epsilon},$$

leads to the term in (53) being

$$\ll \frac{1}{\log K} \sum_{n \leq K^{\delta_0}} \sum_{c \leq K^{1-\delta}} \frac{c^{1+\epsilon}}{c\varphi(c)}\varphi(c) \ll K^{1-\delta/20}. \tag{56}$$

The contribution to (53) from the remainder term for $\chi = \chi_0$ (using (54)) is also at most $K^{1-\delta/2}$ as in (56). We are left with the contribution from main term from $\chi = \chi_0(\mathrm{mod}\, c)$ to (53). That is (53) equals

$$-\frac{3i^a}{\pi \log K} \sum_{m^2 n \leq K^{\delta_0}} \frac{1}{m^2 n} \sum_{c \leq K} \frac{1}{c\varphi(c)} \sum^*_{a(\mathrm{mod}\, c)} S(a, n^2; c)$$

$$\times \int_0^\infty \hat{\phi}\left(\frac{\log t}{2\log K}\right) h\left(\frac{4\pi n\sqrt{t}}{cK}\right) \frac{dt}{\sqrt{t}} + O(K^{1-\delta/2}). \tag{57}$$

Put

$$y = \frac{4\pi n\sqrt{t}}{cK},$$

then (57) is asymptotic to

$$-\frac{3i^a K}{2\pi^2 \log K} \sum_{m^2 n \leq K^{\delta_0}} \frac{1}{m^2 n^2} \sum_{c \leq K} \frac{1}{\varphi(c)} \sum^*_{a(\mathrm{mod}\, c)} S(a, n^2; c) \int_0^\infty \hat{\phi}\left(\frac{\log \frac{Kcy}{4\pi n}}{\log K}\right) h(y) dy.$$

Evaluating the sum of Kloosterman sums, we see the above

$$\sim -\frac{3i^a K\hat{h}(0)}{2\pi^2 \log K} \sum_{m^2 n \leq K^{\delta_0}} \frac{1}{m^2 n^2} \sum_{c \leq K} \frac{\mu(c)}{\varphi(c)} \sum_{d|n^2, d|c} \mu(n^2/d) d\hat{\phi}\left(1 + \frac{\log c}{\log K}\right)$$

$$\sim -\frac{3i^a K\hat{h}(0)}{2\pi^2 \log K} \sum_{m^2 n \leq K^{\delta_0}} \frac{1}{m^2 n^2} \sum_{d|n^2} d \sum_{\lambda \leq K/d} \frac{\mu(\lambda d)}{\varphi(\lambda d)} \mu(\lambda)\hat{\phi}\left(1 + \frac{\log(\lambda d)}{\log K}\right)$$

$$= -\frac{3i^a K\hat{h}(0)}{2\pi^2 \log K} \sum_{m^2 n \leq K^{\delta_0}} \frac{1}{m^2 n^2} \sum_{d|n^2} \frac{d\mu(d)}{\varphi(d)} \sum_{\lambda \leq K/d, (\lambda,d)=1} \frac{\mu^2(\lambda)}{\varphi(\lambda)} \hat{\phi}\left(1 + \frac{\log(\lambda d)}{\log K}\right). \quad (58)$$

Now if

$$F_d(t) = \sum_{\lambda \leq t, (\lambda,d)=1} \frac{\mu^2(\lambda)}{\varphi(\lambda)},$$

then for $d$ square-free

$$F_d(t) \sim \frac{\varphi(d)}{d} \log t, \quad \text{as } t \to \infty.$$

Hence

$$\int_0^\infty \hat{\phi}\left(1 + \frac{\log(dt)}{(\log K)}\right) dF_d(t) \sim \frac{\varphi(d)}{d} \log K \int_0^\infty \hat{\phi}(1+x) dx.$$

Thus (58) is asymptotic to

$$-\frac{3i^a K\hat{h}(0)}{2\pi^2 \log K}(\log K) \sum_{m^2 n \leq K^{\delta_0}} \frac{1}{m^2 n^2} \sum_{d|n^2} \mu(d) \int_0^\infty \hat{\phi}(1+x) dx$$

$$\sim -\frac{3i^a K\hat{h}(0)}{2\pi^2} \frac{\pi^2}{6} \int_0^\infty \hat{\phi}(1+x) dx \sim -\frac{1}{8} i^a K \hat{h}(0) \int_{|\xi| \geqslant 1} \hat{\phi}(\xi) d\xi.$$

This completes the calculation of the contribution of the second term in (50). Collecting all the contributions yields

$$S_K^{(1-a)}[D_K(f, \phi)] \sim \frac{1}{4}\hat{h}(0)K[\int_{-\infty}^\infty \phi(x)dx + \frac{1}{2}\phi(0) - \frac{1}{2}i^a \int_{|\xi| \geqslant 1} \hat{\phi}(\xi) d\xi]. \quad (59)$$

If $a = 0$, this gives

$$S_K^1[D_K(f, \phi)] \sim \frac{1}{4}\hat{h}(0)K[\int_{-\infty}^\infty \phi(x)dx + \frac{1}{2}\int_{|\xi| \leqslant 1} \hat{\phi}(\xi) d\xi]$$

$$= \frac{1}{4}\hat{h}(0)K[\int_{-\infty}^\infty \left(1 + \frac{\sin(2\pi x)}{2\pi x} dx\right) \phi(x) dx]$$

$$\sim S_K^1[1] \int_{-\infty}^\infty \omega(SO(\text{even}))(x)\phi(x)dx. \quad (60)$$

Similarly if $a = 2$, we get

$$S_K^{-1}[D_K(f,\ \phi)] \sim S_K^{-1}[1] \int_{-\infty}^{\infty} \omega(SO(\text{odd}))(x)\phi(x)dx. \tag{61}$$

To deduce Theorem 1 for the family $\mathcal{F}_K$ from (60) and (61) we need to replace $D_K(f,\ \phi)$ by $D(f,\ \phi)$. If we assume, as we are, that the $\gamma_f$'s are real, then for $AK \leqslant k \leqslant BK$ ($A < B$ fixed) and for $\gamma_f > 0$ (say)

$$\phi\left(\gamma_f \frac{\log k}{\pi}\right) = \phi\left(\gamma_f \frac{\log K}{\pi}\right) + O(\gamma_f |\phi'(\gamma_f \xi_f)|), \tag{62}$$

with $\log(AK) \leqslant \xi_f \leqslant \log(BK)$.

Using this and what has been established in (60) and (61) we conclude as in [RS] that in fact

$$S_K^1[D(f,\ \phi)] \sim S_K^1[1] \int_{-\infty}^{\infty} \omega(SO(\text{even}))(x)\phi(x)dx$$

and

$$S_K^{-1}[D(f,\ \phi)] \sim S_K^{-1}[1] \int_{-\infty}^{\infty} \omega(SO(\text{odd}))(x)\phi(x)dx$$

for any $\phi$ with support $(\hat{\phi}) \subset (-2,\ 2)$. This completes the proof of Theorem 1 for the family $\mathcal{F}_K$.

## 2.4 The symmetric square

In this subsection we investigate the distribution of the low lying zeros for the family of $L$-functions $L(s,\ F)$ where $F = \text{sym}^2(f)$, $f \in H_k^*(\Gamma)$. As in the previous sections we let $\rho_F = \frac{1}{2} + i\gamma_F$ denote a typical (non-trivial) zero of $L(s,\ F)$. Throughout this subsection we assume GRH. From (14) and the usual derivation ([RS]) of the explicit formula we obtain, using the relation (9) and with $\phi$ as in Proposition 3 ($AK \leq k \leq BK$) the approximation ($K \to \infty$):

$$\sum_{\gamma_F} \phi\left(\frac{\gamma_F \log K^2}{2\pi}\right) = \int_{-\infty}^{\infty} \phi(x)dx - \frac{1}{\log K}\sum_p \frac{\log p}{p}\hat{\phi}\left(\frac{\log p}{\log K}\right)$$

$$-\frac{1}{\log K} \sum_p \frac{\log p}{\sqrt{p}} \lambda_f(p^2) \hat{\phi}\left(\frac{\log p}{2 \log K}\right)$$

$$-\frac{1}{\log K} \sum_p \frac{\log p}{p} (\lambda_f(p^4) - \lambda_f(p^2)) \hat{\phi}\left(\frac{\log p}{\log K}\right) + O\left(\frac{1}{\log K}\right). \tag{63}$$

Setting

$$D_K(F, \phi) = \phi\left(\frac{\gamma_F \log K}{\pi}\right), \tag{64}$$

we define

$$S_K[D_K(F, \phi)] := \sum_{k \equiv 0 (\bmod 2)} h\left(\frac{k-1}{K}\right) \frac{12}{k-1} \sum_{f \in H_k^*(\Gamma)} D_K(F, \phi) \tag{65}$$

$$= A_K \left[\frac{6}{\pi^2} L(1, F) D_K(F, \phi)\right] \tag{66}$$

The effort of this subsection is directed to proving that for support$(\hat{\phi}) \subset (-4/3, 4/3)$ and $K \to \infty$,

$$S_K[D_K(F, \phi)] \sim \frac{\hat{h}(0) K}{2} \left(\int_\infty^\infty \phi(x) \left(1 - \frac{\sin(2\pi x)}{2\pi x}\right) dx\right)$$

$$\sim S_K[1] \int_\infty^\infty \phi(x) \left(1 - \frac{\sin(2\pi x)}{2\pi x}\right) dx. \tag{67}$$

Now as in (44) we have

$$D_K(F, \phi) \ll \log K. \tag{68}$$

Hence using the approximation (46) together with (66) and (68), we get

$$S_K[D_K(F, \phi)] = \frac{6}{\pi^2} A_K \left[\sum_{m^2 n \leq K^{\delta_0}} \frac{\lambda_f(n^2)}{m^2 n} D_K(F, \phi)\right] + O(K^{1-\delta_0/2}), \tag{69}$$

where $\delta_0 > 0$ is small enough but fixed (depending later on support$(\hat{\phi}) \subset (-4/3, 4/3)$). Using the approximation (63) we have (with the obvious meanings)

$$\frac{6}{\pi^2} A_K \left[\sum_{m^2 n \leq K^{\delta_0}} \frac{\lambda_f(n^2)}{m^2 n} D_K(F, \phi)\right] = I + II + III + IV + V. \tag{70}$$

As in (50) and (51), the sum of the first two terms is asymptotic to

$$I + II \sim \frac{1}{2}\hat{h}(0)K(\int_{\infty}^{\infty} \phi(x)dx - \phi(0)/2), \tag{71}$$

while the last term

$$V = O(K/\log K). \tag{72}$$

The term IV is equal to

$$-\frac{6}{\pi^2 \log K} \sum_p \frac{\log p}{p} \hat{\phi}\left(\frac{\log p}{\log K}\right) A_K \left[\sum_{m^2 n \leq K^{\delta_0}} (\lambda_f(p^4) - \lambda_f(p^2))\frac{\lambda_f(n^2)}{m^2 n}\right]. \tag{73}$$

Applying Propositions 1 and 2 we find that for $\phi$ with support($\hat{\phi}$) $\subset (-5/3, 5/3)$

$$IV \ll K/\log K. \tag{74}$$

So we are left with the term III whose analysis will occupy the rest of this subsection. Remarkably it contributes to the main term (i.e. it is of size $K$) when support($\hat{\phi}$) is not inside $[-1, 1]$. That is this term gives a 'non-diagonal' contribution.

$$-\frac{6}{\pi^2 \log K} \sum_p \frac{\log p}{\sqrt{p}} \hat{\phi}\left(\frac{\log p}{2 \log K}\right) A_K \left[\sum_{m^2 n \leq K^{\delta_0}} \frac{\lambda_f(n^2)\lambda_f(p^2)}{m^2 n}\right] \tag{75}$$

$$= -\frac{6}{\pi^2 \log K} \sum_p \frac{\log p}{\sqrt{p}} \hat{\phi}\left(\frac{\log p}{2 \log K}\right) \sum_{k \equiv 0 (\bmod 2)} h\left(\frac{k-1}{K}\right) \frac{2\pi^2}{k-1}$$

$$\times \sum_{m^2 n \leq K^{\delta_0}} \sum_{f \in H_k^*(\Gamma)} \frac{\lambda_f(n^2)\lambda_f(p^2)}{m^2 n L(1, F)}. \tag{76}$$

By Proposition 2 this equals

$$-\frac{6}{\pi^2 \log K} \sum_p \frac{\log p}{\sqrt{p}} \hat{\phi}\left(\frac{\log p}{2 \log K}\right) \sum_{k \equiv 0 (\bmod 2)} h\left(\frac{k-1}{K}\right)$$

$$\times \sum_{m^2 n \leq K^{\delta_0}} \frac{1}{m^2 n}\left[\delta(p^2, n^2) + 2\pi i^k \sum_{c=1}^{\infty} \frac{S(p^2, n^2; c)}{c} J_{k-1}\left(\frac{4\pi pn}{c}\right)\right] \tag{77}$$

$$= O\left(\frac{K}{\log K}\right) - \frac{12}{\pi \log K} \sum_p \frac{\log p}{\sqrt{p}} \hat{\phi}\left(\frac{\log p}{2 \log K}\right)$$

$$\times \sum_{m^2 n \leq K^{\delta_0}} \frac{1}{m^2 n} \sum_{c=1}^{\infty} \frac{S(p^2, n^2; c)}{c} \sum_{l \equiv 1 (\bmod 2)} i^{l+1} h\left(\frac{l}{K}\right) J_l\left(\frac{4\pi pn}{c}\right). \quad (78)$$

Applying Proposition 1 we have that III is equal to the imaginary part of

$$-\frac{6}{\pi \log K} \sum_p \frac{\log p}{\sqrt{p}} \hat{\phi}\left(\frac{\log p}{2 \log K}\right)$$

$$\times \sum_{m^2 n \leq K^{\delta_0}} \frac{1}{m^2 n} \sum_{c=1}^{\infty} \frac{S(p^2, n^2; c)}{c} e^{i(4\pi pn/c - \pi/4)} \frac{K\sqrt{c}}{\sqrt{4\pi pn}} V_h\left(\frac{cK^2}{8\pi pn}\right) + O\left(\frac{K}{\log K} + \frac{P^{3/2}}{K^{4-\delta_0}}\right), \quad (79)$$

where $P$ is the upper limit on the range of $p$ imposed by the support condition on $\hat{\phi}$.

In particular if

$$\text{support}(\hat{\phi}) \subset (-5/3, \; 5/3), \quad (80)$$

then up to $o(K)$, III is equal to the imaginary part of

$$\frac{-3Ke^{-i\pi/4}}{\pi^{3/2} \log K} \sum_p \frac{\log p}{p} \hat{\phi}\left(\frac{\log p}{2\log K}\right) \sum_{m^2 n \leq K^{\delta_0}} \frac{1}{m^2 n^{3/2}}$$

$$\times \sum_c \frac{S(p^2, n^2; c)}{\sqrt{c}} e\left(\frac{2\pi n}{c}\right) V_h\left(\frac{K^2 c}{8\pi pn}\right) \quad (81)$$

$$= \frac{-3Ke^{-i\pi/4}}{\pi^{3/2} \log K} \sum_{m^2 n \leq K^{\delta_0}} \frac{1}{m^2 n^{3/2}} \sum_c \frac{1}{\sqrt{c}\varphi(c)} \sum_{\chi (\bmod c)}$$

$$\times \sum_{a (\bmod c)}^{*} S(n^2, a^2; c) \bar{\chi}(a) e\left(\frac{2na}{c}\right) \sum_p \frac{\chi(p) \log p}{p} \hat{\phi}\left(\frac{\log p}{2\log K}\right) V_h\left(\frac{cK^2}{8\pi pn}\right)$$

$$= \frac{-3K}{\pi^{3/2} \log K} \sum_{m^2 n \leq K^{\delta_0}} \frac{1}{m^2 n^{3/2}} \sum_c \frac{S_n(c)}{\sqrt{c}\varphi(c)} \sum_p \frac{\log p}{p} \hat{\phi}\left(\frac{\log p}{2\log K}\right) V_h\left(\frac{cK^2}{8\pi pn}\right)$$

$$+ \sum_c \frac{1}{\sqrt{c}\varphi(c)} \sum_{\chi \neq \chi_0} \sum_{a(\bmod c)}^{*} S(n^2, a^2; c) \bar{\chi}(a) e\left(\frac{2na}{c}\right)$$

$$\times \sum_p \frac{\chi(p) \log p}{p} \hat{\phi}\left(\frac{\log p}{2\log K}\right) V_h\left(\frac{cK^2}{8\pi pn}\right),$$

where

$$S_n(c) := \sum_{a(\bmod c)}^{*} S(n^2, a^2; c) e\left(\frac{2na}{c}\right). \quad (82)$$

Now for $\chi \neq \chi_0$ the dyadic sum for $P > K^2$ (the rest in the above being small) is

$$\sum_{p \sim P} \frac{\chi(p) \log p}{p} \hat{\phi}\left(\frac{\log p}{2 \log K}\right) V_h\left(\frac{cK^2}{8\pi pn}\right) \ll_\epsilon P^{-\frac{1}{2}+\epsilon}, \tag{83}$$

using GRH for Dirichlet $L$-functions. Hence the contribution from $\chi \neq \chi_0$ in (83) is

$$\ll K \sum_{m^2 n \leq K^{\delta_0}} \sum_{c \leq Pn/K^2} \frac{c^\epsilon}{\varphi(c)} \frac{c^2}{m^2 n^{3/2}} P^{-\frac{1}{2}+\epsilon} \ll KP^{-\frac{1}{2}+\epsilon} \left(\frac{P}{K^{2-\delta_0}}\right)^{2+\epsilon} = o(K), \tag{84}$$

if support $(\hat{\phi}) \subset (-4/3,\ 4/3)$. Thus we are left with III, which is equal to the imaginary part of

$$\frac{-3Ke^{-i\pi/4}}{\pi^{3/2} \log K} \sum_{m^2 n \leq K^{\delta_0}} \frac{1}{m^2 n^{3/2}} \sum_c \frac{S_n(c)}{\sqrt{c}\varphi(c)} \sum_p \frac{\log p}{p} \hat{\phi}\left(\frac{\log p}{2\log K}\right) V_h\left(\frac{cK^2}{8\pi pn}\right) + o(K) \tag{85}$$

$$= \frac{-3Ke^{-i\pi/4}}{\pi^{3/2} \log K} \sum_{m^2 n \leq K^{\delta_0}} \frac{1}{m^2 n^{3/2}} \sum_c \frac{S_n(c)}{\sqrt{c}\varphi(c)} \sum_{p \geq K^2/n} \frac{\log p}{p} \hat{\phi}\left(\frac{\log p}{2\log K}\right) V_h\left(\frac{cK^2}{8\pi pn}\right) + o(K). \tag{86}$$

The sum $S_n(c)$ may be evaluated directly:

$$S_n(c_1 c_2) = S_n(c_1) S_n(c_2) \quad \text{if } (c_1,\ c_2) = 1; \tag{87}$$

If $(p,\ n) = 1$, then

$$S_n(p^r) = \begin{cases} 1 & \text{if } r = 1 \\ 0 & \text{if } r > 1, 2 \nmid r \\ \varphi(p^r) p^{r/2} & \text{if } r > 1, 2 | r \end{cases} ; \tag{88}$$

If $p|n$, then

$$S_n(p^r) = \begin{cases} -\varphi(p) & \text{if } r = 1 \\ 0 & \text{if } r > 1 \end{cases} . \tag{89}$$

In particular, $S_n(c) = 0$ implies that

$$c = a_1 a_2 b^2, \tag{90}$$

where $p|a_1 \Rightarrow p|n$, $(n, a_2 b^2) = 1$, $(a_2, b) = 1$, $\mu(a_1 a_2) \neq 0$. So the $c$-sum in (87) becomes a sum over $a_1,\ a_2,\ c$, satisfying the above conditions. Now

$$\sum_c \frac{S_n(c)}{\sqrt{c}\varphi(c)} \sum_{p \geq K^2/n} \frac{\log p}{p} \hat{\phi}\left(\frac{\log p}{2\log K}\right) V_h\left(\frac{K^2 c}{8\pi pn}\right)$$

$$= \sum_{a_1, a_2, b} \frac{S_n(a_1)}{\sqrt{a_1}\varphi(a_1)} \frac{S_n(a_2)}{\sqrt{a_2}\varphi(a_2)} \frac{S_n(b^2)}{b\varphi(b^2)} \sum_{p \geq K^2/n} \frac{\log p}{p} \hat{\phi}\left(\frac{\log p}{2\log K}\right) V_h\left(\frac{K^2 a_1 a_2 b^2}{8\pi pn}\right),$$

which according to (89) and (90) is equal to

$$\sum_{a_1,a_2} \frac{\mu(a_1)}{\sqrt{a_1}} \frac{\mu^2(a_2)}{\sqrt{a_2}\varphi(a_2)} \sum_{b\geq 1,(b,na_2)=1} \sum_{p\geq K^2/n} \frac{\log p}{p} \hat{\phi}\left(\frac{\log p}{2\log K}\right) V_h\left(\frac{K^2 a_1 a_2 b^2}{8\pi p n}\right)$$

$$= \sum_{a_1,a_2} \frac{\mu(a_1)}{\sqrt{a_1}} \frac{\mu^2(a_2)}{\sqrt{a_2}\varphi(a_2)} \sum_{r|na_2} \mu(r) \sum_{\nu\geq 1} \sum_{p\geq K^2/n} \frac{\log p}{p} \hat{\phi}\left(\frac{\log p}{2\log K}\right) V_h\left(\frac{K^2 a_1 a_2 r^2 \nu^2}{8\pi p n}\right) \quad (91)$$

Applying Poisson summation to the $\nu$-sum, we deduce that

$$\sum_{\nu\geq 1} V_h\left(\frac{K^2 a_1 a_2 r^2 \nu^2}{8\pi p n}\right) = \frac{1}{2}\left(\sum_{\nu\in\mathbf{Z}} V_h\left(\frac{K^2 a_1 a_2 r^2 \nu^2}{8\pi p n}\right) - V_h(0)\right)$$

$$= \frac{1}{2}\left(\sum_{s\in\mathbf{Z}} \int_{-\infty}^{\infty} V_h\left(\frac{K^2 a_1 a_2 r^2 x^2}{8\pi p n}\right) e(xs) dx - V_h(0)\right). \quad (92)$$

The contribution to (92) of the $-V_h(0)/2$ term is

$$\frac{-V_h(0)}{2} \sum_{a_1,a_2} \frac{\mu(a_1)}{\sqrt{a_1}} \frac{\mu^2(a_2)}{\sqrt{a_2}\varphi(a_2)} \sum_{r|na_2} \mu(r) \sum_{p} \frac{\log p}{p} \hat{\phi}\left(\frac{\log p}{2\log K}\right). \quad (93)$$

Hence its contribution to (87) is

$$\frac{3K V_h(0) e^{-i\pi/4}}{2\pi^{3/2} \log K} \sum_{m^2 n\leq K^{\delta_0}} \frac{1}{m^2 n^{3/2}} \sum_{a_1,a_2} \frac{\mu(a_1)}{\sqrt{a_1}} \frac{\mu^2(a_2)}{\sqrt{a_2}\varphi(a_2)}$$

$$\times \sum_{r|na_2} \mu(r) \sum_{p\geq K^2/n} \frac{\log p}{p} \hat{\phi}\left(\frac{\log p}{2\log K}\right).$$

As $K\to\infty$ we get a contribution to the above only if $na_2 = 1$, i.e. $n=1$ and $a_2 = 1$. So the above is asymptotic to

$$\frac{3K V_h(0) e^{-i\pi/4}}{2\pi^{3/2} \log K} \left(\sum_{m^2\leq K^{\delta_0}} \frac{1}{m^2}\right) \sum_{p\geq K^2/n} \frac{\log p}{p} \hat{\phi}\left(\frac{\log p}{2\log K}\right),$$

which by the prime number theorem is

$$\sim \frac{3K V_h(0) e^{-i\pi/4}}{2\pi^{3/2}} \frac{\pi^2}{6} \int_1^{\infty} \hat{\phi}(y) dy.$$

Recall (see Proposition 1) that $V_h(0) = \frac{\sqrt{2}}{\sqrt{\pi}} \hat{h}(0)$, so the contribution from the $\frac{-V_h(0)}{2}$ term to III is the imaginary part of

$$e^{-i\pi/4} \frac{K}{2} \frac{\sqrt{2}}{\sqrt{\pi}} \sqrt{\pi} \hat{h}(0) \int_1^{\infty} \hat{\phi}(y) dy,$$

which equals
$$\frac{K}{2}\hat{h}(0) - \frac{1}{2}\int_{|y|\geq 1}\hat{\phi}(y)dy. \tag{94}$$

Thus (94) contributes to the main term. We now show that the rest (i.e. the terms in (92) other than $-\frac{1}{2}V_h(0)$) are $o(K)$. We have

$$\sum_{m^2n\leq K^{\delta_0}}\frac{1}{m^2n^{3/2}}\sum_{a_1,a_2}\frac{\mu(a_1)}{\sqrt{a_1}}\frac{\mu^2(a_2)}{\sqrt{a_2}\varphi(a_2)}\sum_{r|na_2}\mu(r)$$
$$\times \sum_{p\geq K^2/n}\frac{\log p}{p}\sum_{s\in\mathbf{Z}}\left(\int_{-\infty}^{\infty}V_h\left(\frac{K^2a_1a_2r^2x^2}{8\pi pn}\right)e(xs)dx\right)\hat{\phi}\left(\frac{\log p}{2\log K}\right). \tag{95}$$

Put $r = r_1r_2$ with $r_1|n$, $r_2|a_2$ and $a_2 = r_2c_2$. Hence (95) becomes

$$\sum_{m^2n\leq K^{\delta_0}}\frac{1}{m^2n^{3/2}}\sum_{r_1\geq 1, r_1|n}\mu(r_1)\sum_{r_2\geq 1,(r_2,n)=1}\frac{\mu(r_2)}{\sqrt{r_2}\varphi(r_2)}\sum_{a_1|n}\frac{\mu(a_1)}{\sqrt{a_1}}\sum_{(c_2,r_2)=1}\frac{\mu(c_2)}{\sqrt{c_2}\varphi(c_2)}$$
$$\times \sum_{p\geq K^2/n}\frac{\log p}{p}\hat{\phi}\left(\frac{\log p}{2\log K}\right)\sum_{s\in\mathbf{Z}}\int_{-\infty}^{\infty}V_h\left(\frac{K^2a_1a_2r^2x^2}{8\pi pn}\right)e(xs)dx. \tag{96}$$

Now according to Proposition 1, for $\beta > 0$

$$\int_{-\infty}^{\infty}V_h(\beta x^2)e(xs)dx = \int_{-\infty}^{\infty}\frac{1}{\sqrt{2\pi}}\int_0^{\infty}\frac{h(\sqrt{u})}{\sqrt{u}}e^{i\beta x^2 u}e(xs)dudx$$
$$= \frac{1+i}{2\sqrt{\beta}}\int_0^{\infty}h(\sqrt{u})e^{-i\pi^2s^2/(\beta u)}\frac{du}{u}. \tag{97}$$

Hence for $|s| \leq \sqrt{\beta}$,

$$\Im\left(e^{-i\pi/4}\int_{-\infty}^{\infty}V_h(\beta x^2)e(xs)dx\right) \ll |s|^2\beta^{-3/2}, \tag{98}$$

and

$$\Im\left(e^{-i\pi/4}\int_{-\infty}^{\infty}V_h(\beta x^2)dx\right) = 0. \tag{99}$$

Also

$$\int_{-\infty}^{\infty}V_h(\beta x^2)e(xs)dx = \frac{1}{\sqrt{\beta}}B\left(\frac{s}{\sqrt{\beta}}\right), \tag{100}$$

where

$$B(v) = \int_{-\infty}^{\infty}V_h(y^2)e(yv)dy. \tag{101}$$

$B(v)$ is rapidly decreasing so from (100) we have that

$$\int_{-\infty}^{\infty} V_h(\beta x^2) e(xs) dx \ll \frac{\sqrt{\beta}}{s^2}. \tag{102}$$

Thus for $\beta > 0$

$$\Im(e^{i\pi/4} \sum_{s \in \mathbf{Z}} \int_{-\infty}^{\infty} V_h(\beta x^2) e(xs) dx)$$
$$= \Im(e^{i\pi/4} \sum_{|s| \leq \beta^{\frac{1}{2}}} \int_{-\infty}^{\infty} V_h(\beta x^2) e(xs) dx) + \Im(e^{i\pi/4} \sum_{|s| > \beta^{\frac{1}{2}}} \int_{-\infty}^{\infty} V_h(\beta x^2) e(xs) dx),$$

which by (99) and (102)

$$\ll \frac{1}{\beta^{3/2}} \sum_{|s| \leq \beta^{\frac{1}{2}}} s^2 + \sum_{|s| > \beta^{\frac{1}{2}}} \frac{\sqrt{\beta}}{s^2} \ll 1. \tag{103}$$

We return to the sum in (96). For $A$ a large parameter we split it into the ranges

$$\frac{K\sqrt{a_1 c_2} r_1 r_2^{3/2}}{\sqrt{8\pi p n}} \leq A$$

and

$$\frac{K\sqrt{a_1 c_2} r_1 r_2^{3/2}}{\sqrt{8\pi p n}} > A$$

In the first case the sum in (96) satisfies (using (102))

$$\Im(e^{-i\pi/4} \sum)$$
$$\ll \sum_{m^2 n \leq K^{\delta_0}} \frac{1}{m^2 n^{3/2}} \sum_{r_1 | n} \sum_{r_2} \frac{1}{\sqrt{r_2} \varphi(r_2)} \sum_{a_1 | n} \frac{1}{\sqrt{a_1}} \sum_{c_2} \frac{1}{\sqrt{c_2} \varphi(c_2)}$$
$$\times \sum_{p \geq K^2 a_1 c_2 r_1^2 r_2^3 / (8\pi n A^2)} \frac{\log p}{p} \sum_{s \neq 0} \frac{K\sqrt{a_1 c_2} r_1 r_2^{3/2}}{s^2 \sqrt{pn}}$$

( $s = 0$ doesn't occur in view of (99) ), which by the prime number theorem

$$\ll \left(\frac{8\pi n A^2}{K^2 a_1 c_2 r_1^2 r_2^3}\right)^{\frac{1}{2}} \frac{K\sqrt{a_1 c_2} r_1 r_2^{3/2}}{\sqrt{n}} \ll A. \tag{104}$$

In the second case we use (103) to conclude that (with an absolute implied constant)

$$\Im(e^{-i\pi/4}\sum) \ll \sum_m \sum_{p \geq \frac{K^2}{n}, r_1|n}{}' \frac{\tau(n)}{m^2 n^{3/2}} \frac{1}{\sqrt{r_2}\varphi(r_2)} \frac{1}{\sqrt{c_2}\varphi(c_2)} \log K, \quad (105)$$

where the stroke $'$ restricts the summation by the additional condition $K\sqrt{a_1 c_2} r_1 r_2^{3/2} \geq A\sqrt{8\pi p n}$. Since $p \geq K^2/n$, it follows that in the sum (105) one of $r_1$ (and hence $n$), $r_2$, $c_2$ or $n$ is large in term of $A$. Since the series involving these all converge it follows that the quantity in (105) is

$$\ll \epsilon_A \log K, \quad (106)$$

where $\epsilon_A \to 0$ as $A \to \infty$. Thus we have shown that

$$\Im(e^{-i\pi/4}\sum) \ll \epsilon_A \log K,$$

and hence combining (106) with (94) we obtain that for support$(\hat{\phi}) \subset (-4/3, 4/3)$

$$III \sim \frac{1}{4} K \hat{h}(0) \int_{|y| \leq 1} \hat{\phi}(y) dy. \quad (107)$$

Therefore the entire contribution to (69) (or (70)) is asymptotic to

$$\frac{\hat{h}(0)}{2} K \left( \int_{-\infty}^{\infty} \phi(x) dx - \phi(0)/2 - \frac{1}{2} \int_{|y| \leq 1} \hat{\phi}(y) dy \right)$$

$$= \frac{\hat{h}(0)}{2} K \int_{-\infty}^{\infty} \left( 1 - \frac{\sin(2\pi x)}{2\pi x} \right) \phi(x) dx. \quad (108)$$

This completes the proof of (67).

Finally from an approximation argument as mentioned in Section 1.3 one derives

$$S_K[D(F, \phi)] \sim S_K[1] \int_{-\infty}^{\infty} \left( 1 - \frac{\sin(2\pi x)}{2\pi x} \right) \phi(x) dx, \quad (109)$$

for support$(\hat{\phi}) \subset (-4/3, 4/3)$.

This concludes the proof of Theorem 1 for the family $\mathcal{F}_{\text{sym}^2}$.

# 3 Large level

## 3.1 Preliminaries

This section is devoted to proving Theorem 1 for the family $\mathcal{F}_N$ of forms of large level $N$ (we assume $N$ is a prime). For simplicity we will examine what is perhaps the most interesting case - that of weight $k = 2$. The analysis goes through in general (i.e. weight $k > 2$, $k$ even), the only modification needed being that in the average over $\mathcal{F}_N$ one also has oldforms (see [AL]) but for $k$ fixed the number of these is $O(1)$ as $N \to \infty$ so that their contribution is insignificant.

For $N$ a large prime let $H^*(N) = H_2^*(N)$ denote the set of $L^2$ normalized weight 2 newforms for $\Gamma_0(N)$ (since $N$ is a prime and there are no such forms for $N = 1$, all Hecke eigenforms of weight 2 for $\Gamma_0(N)$ are newforms). Any $f \in H^*(N)$ has a Fourier expansion

$$f(z) = \sum_{n=1}^{\infty} n^{\frac{1}{2}} a_f(n) e(nz). \tag{110}$$

Moreover for the normalized Hecke operator $T_n$, $n \geqslant 1$,

$$T_n f = \lambda_f(n) f. \tag{111}$$

The Fourier coefficients and the Hecke eigenvalues are related by

$$a_f(n) = a_f(1) \lambda_f(n). \tag{112}$$

Moreover, the eigenvalues $\lambda_f(n)$ satisfy the relations

$$\lambda_f(n) \lambda_f(m) = \sum_{d|(n,m),(d,N)=1} \lambda_f\left(\frac{nm}{d^2}\right). \tag{113}$$

It follows that $f$ is an eigenfunction of the Fricke involution

$$W_n f(z) = N^{-1} z^{-2} f(-1/(Nz)),$$

that is
$$W_n f = \eta_f f, \tag{114}$$

where $\eta_f = \pm 1$, and in terms of the Fourier coefficients

$$\eta_f = \mu(N)\lambda_f(N)N^{\frac{1}{2}}. \tag{115}$$

In terms of the $\lambda_f(n)$'s, the $L$-function $L(s, f)$ is given by

$$L(s, f) = \sum_{n=1}^{\infty} \lambda_f(n) n^{-s} = (1 - \lambda_f(N)N^{-s})^{-1} \prod_{p \neq N} (1 - \lambda_f(p)p^{-s} + p^{-2s})^{-1}. \tag{116}$$

Its functional equation takes the form

$$\Lambda(s, f) = \left(\frac{\sqrt{N}}{2\pi}\right)^s \Gamma(s + \frac{1}{2}) L(s, f) = \epsilon_f \Lambda(1 - s, f), \tag{117}$$

$$\epsilon_f = -\eta_f = -\mu(N)\lambda_f(N)N^{\frac{1}{2}}. \tag{118}$$

In using the Petersson formula (see below) we normalize

$$<f, f>_{\Gamma_0(N)} = \int_{\Gamma_0(N)\backslash \mathbf{H}} |f(z)|^2 y^2 \frac{dxdy}{y^2} = 1. \tag{119}$$

This imposes a normalization on $|a_f(1)|^2$ which we determine in more convenient form (as was done in Section 1). Let $E_\infty(z, s, \Gamma_0(N))$ be the Eisenstein series for the cusp at $\infty$ for $\Gamma_0(N)$,

$$E_\infty(z, s, \Gamma_0(N)) = \sum_{\gamma \in \Gamma_\infty \backslash \Gamma_0(N)} (\Im(\gamma z))^s, \tag{120}$$

where

$$\Gamma_\infty = \left\{ \pm \begin{pmatrix} 1 & t \\ 0 & 1 \end{pmatrix}; \ t \in \mathbf{Z} \right\}.$$

Now for $\Re(s)$ large we get by unfolding method

$$\int_{\Gamma_0(N)\backslash \mathbf{H}} |f(z)|^2 y^2 E_\infty(z, s, \Gamma_0(N)) \frac{dxdy}{y^2}$$
$$= \int_0^\infty \int_0^1 |f(z)|^2 y^{s+1} \frac{dxdy}{y}$$

$$\begin{aligned}
&= \int_0^\infty \sum_{n=1}^\infty n|a_f(n)|^2 e^{-4\pi ny} y^{s+1} \frac{dy}{y} \\
&= |a_f(1)|^2 \sum_{n=1}^\infty \frac{n\lambda_f^2(n)}{(4\pi n)^{s+1}} \Gamma(s+1) \\
&= \frac{|a_f(1)|^2}{(4\pi)^{s+1}} \Gamma(s+1) \sum_{n=1}^\infty \frac{\lambda_f^2(n)}{n^{s+1}} \\
&= \frac{|a_f(1)|^2}{(4\pi)^{s+1}} \Gamma(s+1) \sum_{n=1}^\infty \sum_{d|n,(d,N)=1} \lambda_f\left(\frac{n^2}{d^2}\right) n^{-s} \\
&= \frac{|a_f(1)|^2}{(4\pi)^{s+1}} \Gamma(s+1) \zeta(s)(1-N^{-s}) \sum_{c=1}^\infty \frac{\lambda_f(c^2)}{c^s}.
\end{aligned} \qquad (121)$$

We evaluate the residue at $s=1$ of both sides of (121). First from the general principles ([Sa])

$$\operatorname{Res}_{s=1} E_\infty(z,\ s,\ \Gamma_0(N)) = \frac{1}{\operatorname{Vol}(\Gamma_0(N)\backslash \mathbf{H})} = \frac{3}{\pi(N+1)}. \qquad (122)$$

Imposing (120) we have

$$\frac{3}{\pi(N+1)} = \frac{|a_f(1)|^2}{(4\pi)^2}(1-N^{-1}) \sum_{n=1}^\infty \frac{\lambda_f(n^2)}{n}. \qquad (123)$$

Thus

$$|a_f(1)|^2 = \frac{48\pi}{N-N^{-1}} \left(\sum_{n=1}^\infty \frac{\lambda_f(n^2)}{n}\right)^{-1}. \qquad (124)$$

We recall the Petersson formula (with $<f,\ f>_{\Gamma_0(N)} = 1$).

**Proposition 1** ([Pet]). For $m,n \geqslant 1$, we have

$$(4\pi)^{-1} \sum_{f \in H^*(N)} a_f(m)\overline{a_f(n)} = \delta(m,\ n) - 2\pi \sum_{c \equiv 0(\bmod N)} \frac{S(m,n;c)}{c} J_1\left(\frac{4\pi\sqrt{mn}}{c}\right). \qquad (125)$$

In terms of the $\lambda_f(n)$'s, this becomes (using (124))

$$\frac{12}{N-N^{-1}} \sum_{f \in H^*(N)} \omega_f \lambda_f(n)\lambda_f(m) = \delta(m,\ n) - 2\pi \sum_{c \equiv 0(\bmod N)} \frac{S(m,n;c)}{c} J_1\left(\frac{4\pi\sqrt{mn}}{c}\right), \qquad (126)$$

where
$$\omega_f = (\sum_{t=1}^{\infty} \frac{\lambda_f(t^2)}{t})^{-1}. \tag{127}$$

We define the averaging operator $A^*$ on $H^*(N)$ by

$$A^*[X_f] = \sum_{f \in H^*(N)} \omega_f X_f. \tag{128}$$

As in Section 1, $\omega_f > 0$, and assuming RH for $L(s, \text{sym}^2(f))$ (which we do) we have that for $\delta_0 > 0$ (fixed but arbitrarily small)

$$\omega_f^{-1} = \sum_{t \leq N^{\delta_0}} \frac{\lambda_f(t^2)}{t} + O_\epsilon(N^{-\delta_0/2+\epsilon}). \tag{129}$$

Then formula (126) becomes

$$\frac{12}{N - N^{-1}} A^*[\lambda_f(n)\lambda_f(m)] = \delta(m, n) - 2\pi \sum_{c \equiv 0(\text{ mod } N)} \frac{S(m, n; c)}{c} J_1\left(\frac{4\pi\sqrt{mn}}{c}\right). \tag{130}$$

Applying Weil's bound (26) for Kloosterman sums and the bound $J_1(x) \ll x$, we obtain

**Proposition 2.** For $m, n \geq 1$ we have

$$A^*[\lambda_f(n)\lambda_f(m)] = \frac{N}{12}\delta(m, n) + O_\epsilon(N^\epsilon (m, n)^{\frac{1}{2}}(mn/N)^{\frac{1}{2}}). \tag{131}$$

Factoring the Kloosterman sum $S(N, n; Nl)$ as a Ramanujan sum to the modulus $N$ and a Kloosterman sum to the modulus $l$ (see below), we get in a similar way

**Proposition 3.** We have

$$A^*[\lambda_f(N)\lambda_f(n)] = \frac{N}{12}\delta(N, n) + O_\epsilon(N^\epsilon (n, N)(n/N)^{\frac{1}{2}}). \tag{132}$$

Set

$$S[X_f] = \sum_{f \in H^*(N)} X_f, \tag{133}$$

$$S^+[X_f] = \sum_{f \in H^*(N)} \frac{1 + \epsilon_f}{2} X_f, \tag{134}$$

$$S^-[X_f] = \sum_{f \in H^*(N)} \frac{1-\epsilon_f}{2} X_f. \tag{135}$$

Thus

$$S[X_f] = A^*[X_f \sum_{n=1}^{\infty} \frac{\lambda_f(n^2)}{n}], \tag{136}$$

which if $X_f \ll \log N$, say, satisfies

$$S[X_f] = A^*[X_f \sum_{n \leq N^{\delta_0}} \frac{\lambda_f(n^2)}{n}] + O(N^{1-\delta_0/4}) \tag{137}$$

(using (130) and where $\delta_0 > 0$ is small but fixed). In particular, applying Proposition 2 we have

$$S[1] = \frac{N}{12} + O(N^{1-\delta_0/4}), \tag{138}$$

and applying Proposition 3 and (118)

$$S^{\pm}[1] = \frac{N}{24} + O(N^{1-\delta_0/4}). \tag{139}$$

Finally we recall the approximate explicit formula ([RS]). As in the rest of this paper we write the zeros of $\Lambda(s, f)$ as $\frac{1}{2} + i\gamma_f$.

**Proposition 4.** Let $\phi$ be a test function as in Proposition 3 (of Section 1), then we have

$$\begin{aligned} D(f, \phi) &:= \sum_{\gamma_f} \phi\left(\frac{\log N}{2\pi} \gamma_f\right) \\ &= \int_{-\infty}^{\infty} \phi(x) dx + \frac{2}{\log N} \sum_{p \neq N} \frac{\log p}{p} \hat{\phi}\left(\frac{2 \log p}{\log N}\right) \\ &\quad - \frac{2}{\log N} \sum_{p \neq N} \frac{\lambda_f(p^2) \log p}{p} \hat{\phi}\left(\frac{2 \log p}{\log N}\right) \\ &\quad - \frac{2}{\log N} \sum_{p \neq N} \frac{\lambda_f(p) \log p}{\sqrt{p}} \hat{\phi}\left(\frac{\log p}{\log N}\right) + O\left(\frac{1}{\log N}\right). \end{aligned}$$

## 3.2 Proof of Theorem 1 for $\mathcal{F}_N$

In what follows we examine the asymptotics of

$$S^{\pm}[D(f, \phi)], \tag{140}$$

as $N \to \infty$. We assume that support($\hat{\phi}$) $\subset (-2, 2)$. In particular support($\hat{\phi}$) $\subset (-2+4\delta_0, 2-4\delta_0)$, say, for a suitable $\delta_0 > 0$. We are also assuming GRH throughout (for $L(s, \text{sym}^2(f))$ this was already assumed in (129), while we will assume it for $L(s, \chi)$, $\chi$ a Dirichlet character, later in this Section). According to Proposition 4 we have

$$S^+[D(f, \phi)] = I + II + III + IV + V, \tag{141}$$

with the obvious interpretations of these terms. From (139) we have

$$I + II + V \sim \frac{N}{24} \left( \int_{-\infty}^{\infty} \phi(x) dx + \frac{2}{\log N} \sum_{p \neq N} \frac{\log p}{p} \hat{\phi}\left(\frac{2 \log p}{\log N}\right) \right) + O\left(\frac{N}{\log N}\right). \tag{142}$$

Applying the prime number theorem to the sum above one gets

$$I + II + V \sim \frac{N}{24} \left( \int_{-\infty}^{\infty} \phi(x) dx + \frac{\phi(0)}{2} \right). \tag{143}$$

Using Propositions 2 and 3 and (118) we can easily deal with III as follows.

$$\begin{aligned}
III &\sim -\frac{2}{\log N} \sum_{p \neq N} \frac{\log p}{p} \hat{\phi}\left(\frac{2 \log p}{\log N}\right) A^*[\sum_{n \leq N^{\delta_0}} \frac{\lambda_f(n^2)}{n} \left(\frac{1+\epsilon_f}{2}\right) \lambda_f(p^2)] \\
&= -\frac{2}{\log N} \sum_{p \neq N, n \leq N^{\delta_0}} \frac{\log p}{pn} \hat{\phi}\left(\frac{2 \log p}{\log N}\right) \\
&\quad \times \left[ \frac{N}{24} \delta(n^2, p^2) + O\left((n^2, p^2)^{\frac{1}{2}} \frac{np}{\sqrt{N}}\right) + O_\epsilon\left(N^{\frac{1}{2}} \left(\frac{n^2 p^2}{N}\right)^{\frac{1}{2}+\epsilon}\right) \right] \\
&= -\frac{2N}{24 \log N} \sum_{p \leq N^{\delta_0}} \frac{\log p}{p^2} \hat{\phi}\left(\frac{2 \log p}{\log N}\right) + O(\sum_{p \leq N^{1-2\delta_0}} \log p).
\end{aligned}$$

(since support($\hat{\phi}$) $\subset [-2+4\delta_0, 2-4\delta_0]$). Hence

$$III = o(N). \tag{144}$$

We are left with dealing with the term IV which in fact contributes to the main term (i.e. is of size $N$). This term is the 'non-diagonal' term mentioned in the introduction and which is present if support($\hat{\phi}$) is not contained in $[-1, 1]$.

Now

$$IV = -\frac{2}{\log N} \sum_{p \neq N} \sum_{n \leq N^{\delta_0}} \frac{\log p}{\sqrt{p}n} \hat{\phi}\left(\frac{\log p}{\log N}\right) A^*[\frac{1+\epsilon_f}{2}\lambda_f(n^2)\lambda_f(p)]. \tag{145}$$

We split this as $IV_A + IV_B$ where

$$IV_A = -\frac{1}{\log N} \sum_{p \neq N} \sum_{n \leq N^{\delta_0}} \frac{\log p}{\sqrt{p}n} \hat{\phi}\left(\frac{\log p}{\log N}\right) A^*[\lambda_f(n^2)\lambda_f(p)], \tag{146}$$

and $IV_B$ is the term involving $\epsilon_f$.

Apply (130) to get

$$\begin{aligned}IV_A &= (1+O(N^{-1}))\frac{\pi N}{6\log N} \sum_{p \neq N} \sum_{n \leq N^{\delta_0}} \frac{\log p}{\sqrt{p}n} \hat{\phi}\left(\frac{\log p}{\log N}\right) \\ &\quad \times \sum_{c \equiv 0(\bmod N)} \frac{S(n^2, p; c)}{c} J_1\left(\frac{4\pi\sqrt{p}n}{c}\right).\end{aligned} \tag{147}$$

Up to a lower order term (using $p \leq N^{2-4\delta_0}$ following from support($\hat{\phi}$) $\subset (-2+4\delta_0, 2-4\delta_0)$, and $J_1(z) = z + O(|z|^2)$), we get

$$\begin{aligned}IV_A &= \frac{\pi N}{6\log N} \sum_{p \neq N} \sum_{n \leq N^{\delta_0}} \frac{\log p}{\sqrt{p}n} \hat{\phi}\left(\frac{\log p}{\log N}\right) \sum_{l=1}^{\infty} \frac{S(n^2, p; lN)}{lN} \frac{4\pi n\sqrt{p}}{lN} \\ &= \frac{2\pi^2}{3N\log N} \sum_{p \neq N} \sum_{n \leq N^{\delta_0}} (\log p)\, \hat{\phi}\left(\frac{\log p}{\log N}\right) \sum_{l=1}^{\infty} \frac{S(n^2, p; lN)}{l^2} \\ &= \frac{2\pi^2}{3N\log N} \sum_{p \neq N} \sum_{n \leq N^{\delta_0}} \sum_{l \leq l_0} (\log p)\, \hat{\phi}\left(\frac{\log p}{\log N}\right) \frac{S(n^2, p; lN)}{l^2} \\ &\quad + O_\epsilon\left(N^{-1} \sum_{l \geq l_0} \sum_{p \leq N^2} (\log p) \sum_{n \leq N^{\delta_0}} \frac{(lN)^{\frac{1}{2}+\epsilon}}{l^2}\right),\end{aligned} \tag{148}$$

with $l_0$ to be chosen momentarily.

The above is

$$\frac{2\pi^2}{3N\log N} \sum_{n \leq N^{\delta_0}} \sum_{l \leq l_0} {\sum_{a(\bmod lN)}}^* \frac{S(n^2, a; lN)}{l^2}$$

$$\times \sum_{p \equiv a(\bmod lN)} (\log p)\, \hat{\phi}\left(\frac{\log p}{\log N}\right) + O(N^2 l_0^{-\frac{1}{2}}). \tag{149}$$

We choose $l_0 = N^4$ so that the $O$-term is insignificant. Thus (up to the negligible error term)

$$IV_A = \frac{2\pi^2}{3N \log N} \sum_{n \leq N^{\delta_0}} \sum_{l \leq l_0} l^{-2} \sum_{a(\bmod lN)}^* S(n^2, a; lN)$$

$$\times \frac{1}{\varphi(lN)} \sum_{\chi(\bmod lN)} \bar{\chi}(a) \sum_p \chi(p)(\log p) \,\hat{\phi}\!\left(\frac{\log p}{\log N}\right). \qquad (150)$$

For $\chi \neq \chi_0 (\bmod lN)$ in (150) we assume RH so that

$$\sum_p \chi(p)(\log p) \,\hat{\phi}\!\left(\frac{\log p}{\log N}\right) \ll_\epsilon N^{1-2\delta_0+\epsilon} \qquad (151)$$

(using $p \ll N^{2-4\delta_0}$). Thus the contribution to $IV_A$ from the terms with $\chi \neq \chi_0$ is

$$\ll N^{-1} \sum_{n \leq N^{\delta_0}} \sum_{l \leq l_0} \frac{1}{l^2 \varphi(lN)} \sum_{\chi \neq \chi_0(\bmod lN)} N^{1-2\delta_0+\epsilon} |\sum_{a(\bmod lN)}^* S(n^2, a; lN) \bar{\chi}(a)|$$

$$\ll \frac{N^{1-2\delta_0+\epsilon}}{N} \sum_{n \leq N^{\delta_0}} \sum_{l \leq l_0} \frac{lN}{l^2} \ll N^{1-\delta_0+\epsilon} l_0^\epsilon \ll N^{1-\delta_0+5\epsilon}. \qquad (152)$$

That is the contribution from the $\chi \neq \chi_0$ terms to $IV_A$ is $o(N)$. The contribution to $IV_A$ from $\chi = \chi_0$ is easily estimated as being $o(N)$ as well. Therefore we have

$$IV_A = o(N). \qquad (153)$$

We turn to $IV_B$ which is given as

$$IV_B = \frac{-1}{\log N} \sum_{p \neq N} \sum_{n \leq N^{\delta_0}} \frac{\log p}{\sqrt{p}n} \hat{\phi}\!\left(\frac{\log p}{\log N}\right) A^*[-\mu(N) N^{\frac{1}{2}} \lambda_f(N) \lambda_f(n^2) \lambda_f(p)] \qquad (154)$$

$$= \frac{\mu(N) N^{\frac{1}{2}}}{\log N} \sum_{p \neq N} \sum_{n \leq N^{\delta_0}} \frac{\log p}{\sqrt{p}n} \hat{\phi}\!\left(\frac{\log p}{\log N}\right) A^*[\lambda_f(Nn^2) \lambda_f(p)]$$

$$\sim \frac{N}{12} \frac{\mu(N) N^{\frac{1}{2}}}{\log N} \sum_{p \neq N} \sum_{n \leq N^{\delta_0}} \frac{\log p}{\sqrt{p}n} \hat{\phi}\!\left(\frac{\log p}{\log N}\right)$$

$$\times \left[ -2\pi \sum_{c \equiv 0(\bmod N)}^\infty \frac{S(Nn^2, p; c)}{c} J_1\!\left(\frac{4\pi n \sqrt{Np}}{c}\right) \right]$$

$$= \frac{-\pi\mu(N)N^{\frac{1}{2}}}{6\log N} \sum_{p\neq N} \sum_{n\leq N^{\delta_0}} \frac{\log p}{\sqrt{p}n} \hat{\phi}\left(\frac{\log p}{\log N}\right) \sum_{l=1}^{\infty} \frac{S(Nn^2, p; lN)}{l} J_1\left(\frac{4\pi n\sqrt{p}}{l\sqrt{N}}\right). \qquad (155)$$

The total contribution of the terms in (155) coming from $l$'s with $(l, N) > 1$ or $l \geqslant N^3$, is easily estimated as $o(N)$. Thus up to $o(N)$ we find that $IV_B$ is asymptotically

$$\frac{-\pi\mu(N)N^{\frac{1}{2}}}{6\log N} \sum_{p\neq N} \sum_{n\leq N^{\delta_0}} \frac{\log p}{\sqrt{p}n} \hat{\phi}\left(\frac{\log p}{\log N}\right) \sum_{(l,N)=1,l\leq N^3} \frac{S(Nn^2, p; lN)}{l} J_1\left(\frac{4\pi n\sqrt{p}}{l\sqrt{N}}\right). \qquad (156)$$

Now for $(l, N) = 1$,

$$S(Nn^2, p; lN) = \mu(N)S(n^2, \bar{N}p; l), \qquad (157)$$

where $N\bar{N} \equiv 1 (\mathrm{mod}\, l)$, so

$$\begin{aligned}
IV_B &= \frac{-\pi N^{\frac{1}{2}}}{6\log N} \sum_{p\neq N} \sum_{n\leq N^{\delta_0}} \frac{\log p}{\sqrt{p}n} \hat{\phi}\left(\frac{\log p}{\log N}\right) \sum_{(l,N)=1,l\leqslant N^3} \frac{S(n^2, \bar{N}p; l)}{l} J_1\left(4\pi \frac{n\sqrt{p}}{l\sqrt{N}}\right) \\
&\sim \frac{-\pi N^{\frac{1}{2}}}{6\log N} \sum_{(l,N)=1,l\leqslant N^3} \frac{1}{l} \sum_{n\leq N^{\delta_0}} \frac{1}{n} \sum\nolimits^*_{a(\mathrm{mod}\, l)} S(n^2, \bar{N}a; l) \\
&\quad \times \sum_{p\equiv a(l)} \frac{\log p}{\sqrt{p}} \hat{\phi}\left(\frac{\log p}{\log N}\right) J_1\left(\frac{4\pi n\sqrt{p}}{l\sqrt{N}}\right) \\
&= \frac{-\pi N^{\frac{1}{2}}}{6\log N} \sum_{(l,N)=1,l\leqslant N^3} \sum_{n\leq N^{\delta_0}} \frac{1}{nl\varphi(l)} \sum_{\chi(\mathrm{mod}\, l)} \left(\bar{\chi}(a) \sum\nolimits^*_{a(\mathrm{mod}\, l)} S(n^2, \bar{N}a; l)\right) \\
&\quad \times \sum_{p} \frac{\chi(p)\log p}{\sqrt{p}} \hat{\phi}\left(\frac{\log p}{\log N}\right) J_1\left(\frac{4\pi n\sqrt{p}}{l\sqrt{N}}\right). \qquad (158)
\end{aligned}$$

We separate the terms with $\chi = \chi_0(\mathrm{mod}\, l)$ and those with $\chi \neq \chi_0(\mathrm{mod}\, l)$. In the latter case we have, assuming RH for $L(s, \chi)$

$$\sum_p \frac{\chi(p)\log p}{\sqrt{p}} \hat{\phi}\left(\frac{\log p}{\log N}\right) J_1\left(\frac{4\pi n\sqrt{p}}{l\sqrt{N}}\right) \ll \frac{nN^{1-2\delta_0+\epsilon}}{\sqrt{N}l}. \qquad (159)$$

Hence this contribution is

$$\ll N^{\frac{1}{2}} \sum_{n\leqslant N^{\delta_0}} \frac{1}{n} \sum_{l\leq N^3} \frac{1}{l\varphi(l)} \sum_{\chi\neq\chi_0(\mathrm{mod}\, l)} l^{1+\epsilon} \frac{nN^{1-2\delta_0+\epsilon}}{\sqrt{N}l} \ll N^{1-\delta_0+\epsilon},$$

which is $o(N)$ as $N \to \infty$. Thus up to $o(N)$ we are left with

$$\begin{aligned}
IV_B &= \frac{-\pi N^{\frac{1}{2}}}{6 \log N} \sum_{(l,N)=1, l \leq N^3} \sum_{n \leq N^{\delta_0}} \frac{1}{nl\varphi(l)} \sum^*_{a(\bmod l)} S(n^2, \bar{N}a; l) \\
&\quad \times \sum_p \frac{\log p}{\sqrt{p}} \hat{\phi}\left(\frac{\log p}{\log N}\right) J_1\left(\frac{4\pi n \sqrt{p}}{l\sqrt{N}}\right) \\
&= \frac{-\pi N^{\frac{1}{2}}}{6 \log N} \sum_{(l,N)=1, l \leqslant N^3} \sum_{n \leq N^{\delta_0}} \frac{\mu(l)}{nl\varphi(l)} \sum_{d|n^2, d|l^2} \mu\left(\frac{l}{d}\right) d \\
&\quad \times \sum_p \frac{\log p}{\sqrt{p}} \hat{\phi}\left(\frac{\log p}{\log N}\right) J_1\left(\frac{4\pi n \sqrt{p}}{l\sqrt{N}}\right) \\
&= \frac{-\pi N^{\frac{1}{2}}}{6 \log N} \sum_{n \leqslant N^{\delta_0}} \frac{1}{n} \sum_{d|n^2} \sum_{(dc,N)=1, dc \leqslant N^3} \frac{\mu(dc)d\mu(c)}{dc\varphi(dc)} \\
&\quad \times \sum_p \frac{\log p}{\sqrt{p}} \hat{\phi}\left(\frac{\log p}{\log N}\right) J_1\left(\frac{n\sqrt{p}4\pi}{dc\sqrt{N}}\right). \quad (160)
\end{aligned}$$

By the prime number theorem with the square root remainder we have that up to $o(N)$ the last expression is

$$\sim \frac{-\pi N^{\frac{1}{2}}}{6 \log N} \sum_{n \leqslant N^{\delta_0}} \frac{1}{n} \sum_{d|n^2} \sum_{(dc,N)=1, dc \leqslant N^3} \frac{\mu(dc)\mu(c)}{c\varphi(dc)}$$
$$\times \int_1^\infty J_1\left(\frac{4\pi n \sqrt{t}}{dc\sqrt{N}}\right) \hat{\phi}\left(\frac{\log t}{\log N}\right) \frac{dt}{\sqrt{t}}.$$

Set
$$y = \frac{4\pi n \sqrt{t}}{dc\sqrt{N}},$$

then the above equals

$$\frac{-N}{24 \log N} \sum_{n \leq N^{\delta_0}} \frac{1}{n^2} \sum_{d|n^2} \sum_{(dc,N)=1, dc \leq N^3} \frac{\mu(dc)\mu(c)dc}{c\varphi(dc)} \int_1^\infty J_1(y) \hat{\phi}\left(\frac{2\log(y\sqrt{N}cd/4\pi n)}{\log N}\right) dy. \quad (161)$$

Now the $n$-series converges absolutely so we consider first

$$\sum_{(c,d)=1} \frac{\mu^2(c)}{\varphi(c)} \int_0^\infty J_1(y) \hat{\phi}\left(1 + \frac{\log(ydc/4\pi n)}{\log N}\right) dy$$
$$\sim \sum_{(c,d)=1} \frac{\mu^2(c)}{\varphi(c)} \int_0^\infty J_1(y) \hat{\phi}\left(1 + \frac{\log(yc)}{\log N}\right) dy$$

$$= \int_0^\infty J_1(y) \left( \sum_{(c,d)=1} \frac{\mu^2(c)}{\varphi(c)} \hat{\phi}\left(1 + \frac{\log(yc)}{\log N}\right) \right) dy. \tag{162}$$

As in the discussion following (58) the term in parenthesis is

$$\sim \frac{\varphi(d)}{d} \log N \int_0^\infty \hat{\phi}(1+x) dx.$$

So the above (162) is

$$\sim \frac{\varphi(d)}{d}(\log N) \int_0^\infty J_1(y) \int_0^\infty \hat{\phi}(1+x) dx dy$$
$$= \frac{\varphi(d)}{d}(\log N) \int_1^\infty \hat{\phi}(\xi) d\xi,$$

the $y$ integral being equal to 1 (see (6.511.1) of [GR]). Hence (160) is asymptotic to

$$-\frac{N}{24} \sum_{n \leq N^{\delta_0}} \frac{1}{n^2} \sum_{d|n^2} \frac{d\mu(d)}{\varphi(d)} \frac{\varphi(d)}{d} \int_1^\infty \hat{\phi}(\xi) d\xi \sim -\frac{N}{48} \int_{|\xi| \geq 1} \hat{\phi}(\xi) d\xi.$$

We have shown that $IV_B$ and hence $IV$ satisfies

$$IV \sim -\frac{N}{48} \int_{|\xi| \geq 1} \hat{\phi}(\xi) d\xi. \tag{163}$$

Combining all the terms in (141) we get

$$S^+[D(f, \phi)] \sim \frac{N}{24} \left( \int_{-\infty}^\infty \phi(x) + \frac{\phi(0)}{2} - \frac{1}{2} \int_{|\xi| \geq 1} \hat{\phi}(\xi) d\xi \right)$$
$$= \frac{N}{24} \left( \int_{-\infty}^\infty \phi(x) + \frac{1}{2} \int_{|\xi| \leq 1} \hat{\phi}(\xi) d\xi \right)$$
$$= \frac{N}{24} \int_{-\infty}^\infty \phi(x) \left( 1 + \frac{\sin(2\pi x)}{2\pi x} \right) dx$$
$$\sim S^+[1] \int_{-\infty}^\infty \phi(x) W(SO(\text{even}))(x) dx. \tag{164}$$

For $S^-[D(f, \phi)]$ the only change is that

$$-\frac{1}{2} \int_{|\xi| \geq 1} \hat{\phi}(\xi) d\xi$$

appears as

$$\frac{1}{2} \int_{|\xi| \geq 1} \hat{\phi}(\xi) d\xi.$$

Thus

$$S^-[D(f, \phi)] \sim \frac{N}{24}\left(\int_{-\infty}^{\infty} \phi(x) + \phi(0) - \frac{1}{2}\int_{|\xi|\leq 1} \hat{\phi}(\xi)d\xi\right)$$

$$= \frac{N}{24}\left(\int_{-\infty}^{\infty} \phi(x)\left(1 - \frac{\sin(2\pi x)}{2\pi x}\right)dx + \phi(0)\right)$$

$$\sim S^-[1]\int_{-\infty}^{\infty} \phi(x)W(SO(\mathrm{odd}))(x)dx. \tag{165}$$

The results (164) and (165) are proven under GRH and the assumption that support$(\hat{\phi}) \subset (-2, 2)$. This concludes the proof of Theorem 1 in the $N$-aspect.

# 4 Application to nonvanishing at $s = \frac{1}{2}$

In this short Section we derive Corollary 2 from Theorem 1. To this end define $\alpha_2(G)$, where $G$ is one of $SO(\mathrm{even})$, $SO(\mathrm{odd})$ or $O$, and $W(G)(x)$ is the corresponding density,

$$\alpha_2(G) = \inf \int_{-\infty}^{\infty} W(G)(x)\phi(x)dx, \tag{166}$$

the infimum being over all $\phi \in \mathcal{S}(\mathbf{R})$ with $\phi(x) \geq 0$, $\phi(0) = 1$ and support$(\hat{\phi}) \subset (-2, 2)$. In the Appendix we determine these numbers $\alpha_2(G)$. The key point that emerges is that as long as

$$\widehat{W(G)}(\xi) = \delta_0 + m_G(\xi),$$

where $m_G(\xi)$ is not constant on $[-2, 2]$ (in our example this results from the non-diagonal analysis carried out in the previous Sections), the familiar test function

$$\phi(x) = \left(\frac{\sin 2\pi x}{2\pi x}\right)^2$$

is not optimal. That is, it does not achieve the infimum in (166).

Applying Theorem 1 (assuming GRH as is done in that Theorem) with a test function $\phi$ as in (1) we obtain:

$$\varlimsup_{K\to\infty} \frac{1}{M^+(K)} \sum_{k\equiv 0(\bmod 4), k\leq K} \sum_{f\in H_k(\Gamma)} \mathrm{ord}(\frac{1}{2}, f) \leq \alpha_2(SO(\mathrm{even})); \tag{167}$$

$$\varlimsup_{K \to \infty} \frac{1}{M^-(K)} \sum_{k \equiv 2(\bmod\ 4), k \leqslant K} \sum_{f \in H_k(\Gamma)} \operatorname{ord}(\tfrac{1}{2}, f) \leqslant \alpha_2(SO(\text{odd})); \tag{168}$$

$$\varlimsup_{N \to \infty} \frac{1}{M^+(N)} \sum_{f \in H_k^+(N)} \operatorname{ord}(\tfrac{1}{2}, f) \leqslant \alpha_2(SO(\text{even})); \tag{169}$$

$$\varlimsup_{N \to \infty} \frac{1}{M^-(N)} \sum_{f \in H_k^-(N)} \operatorname{ord}(\tfrac{1}{2}, f) \leqslant \alpha_2(SO(\text{odd})). \tag{170}$$

It is easy to see that in determining the infimum in (1) we can, without loss of generality, allow functions $\phi \in L^1(\mathbf{R})$ with support($\hat{\phi}$) $\subset [-2, 2]$. In particular

$$\phi(x) = \left(\frac{\sin 2\pi x}{2\pi x}\right)^2$$

may be used. However since it does not achieve the infimum we have the strict inequalities

$$\alpha_2(SO(\text{even})) < \int_{-\infty}^{\infty} \left(1 + \frac{\sin 2\pi x}{2\pi x}\right) \left(\frac{\sin 2\pi x}{2\pi x}\right)^2 dx = \frac{7}{8}. \tag{171}$$

$$\alpha_2(SO(\text{odd})) < \int_{-\infty}^{\infty} \left(\delta_0 + \left(1 - \frac{\sin 2\pi x}{2\pi x}\right)\right) \left(\frac{\sin 2\pi x}{2\pi x}\right)^2 dx = \frac{9}{8}. \tag{172}$$

On the other hand $m_O(\xi)$ is constant so that

$$\alpha_2(O) = \int_{-\infty}^{\infty} \left(\frac{\delta_0}{2} + dx\right) \left(\frac{\sin 2\pi x}{2\pi x}\right)^2 dx = 1. \tag{173}$$

From (169) and (171) we have that

$$\varlimsup_{N \to \infty} \frac{1}{M^+(N)} \sum_{\substack{f \in H_k^+(N) \\ \operatorname{ord}(\frac{1}{2}, f) \geqslant 2}} 2 < \frac{7}{8}, \tag{174}$$

and hence

$$\varlimsup_{N \to \infty} \frac{1}{M^+(N)} \#\{f \in H_k^+(N) | \operatorname{ord}(\tfrac{1}{2}, f) \geqslant 2\} < \frac{7}{16}. \tag{175}$$

Since $\operatorname{ord}(\tfrac{1}{2}, f)$ is even for $f \in H_k^+(N)$ we conclude from (175) that

$$\varliminf_{N \to \infty} \frac{1}{M^+(N)} \#\{f \in H_k^+(N) | \operatorname{ord}(\tfrac{1}{2}, f) \neq 0\} > \frac{9}{16}. \tag{176}$$

This establishes part (c) of Corollary 2.

For $M_k^-(N)$ we have from (170) and (172) that

$$\overline{\lim_{N\to\infty}}\frac{1}{M^-(N)}\left[\sum_{f\in H_k^-(N),\,\mathrm{ord}(\frac{1}{2},\,f)=1} 1 + \sum_{f\in H_k^-(N),\,\mathrm{ord}(\frac{1}{2},\,f)\geq 3} 3\right] < \frac{9}{8}. \tag{177}$$

Hence,

$$\overline{\lim_{N\to\infty}}\frac{1}{M^-(N)} \sum_{\substack{f\in H_k^-(N) \\ \mathrm{ord}(\frac{1}{2},\,f)\geq 3}} 2 < \frac{1}{8}, \tag{178}$$

or

$$\overline{\lim_{N\to\infty}}\frac{1}{M^-(N)} \sum_{\substack{f\in H_k^-(N) \\ \mathrm{ord}(\frac{1}{2},\,f)\geq 3}} 1 < \frac{1}{16}. \tag{179}$$

Therefore it follows that

$$\underline{\lim}_{N\to\infty}\frac{1}{M^-(N)}\#\{f\in H_k^-(N)\,|\,\mathrm{ord}(\frac{1}{2},\,f)=1\} > \frac{15}{16}. \tag{180}$$

This proves part (d) of Corollary (2). Parts (a) and (b) are proven in the same way.

Concerning part (c) we have shown above that

$$\overline{\lim_{N\to\infty}}\frac{1}{M^+(N)} \sum_{f\in H_k^+(N)} \mathrm{ord}(\frac{1}{2},\,f) \leq \frac{7}{8}, \tag{181}$$

and

$$\overline{\lim_{N\to\infty}}\frac{1}{M^-(N)} \sum_{f\in H_k^-(N)} \mathrm{ord}(\frac{1}{2},\,f) \leq \frac{9}{8}. \tag{182}$$

It follows on adding these and from

$$M_k^-(N) \sim M_k^+(N) \sim \frac{M_k^*(N)}{2}$$

that

$$\overline{\lim_{N\to\infty}}\frac{1}{M^*(N)} \sum_{f\in H_k^*(N)} \mathrm{ord}(\frac{1}{2},\,f) < \frac{1}{2}\left(\frac{7}{8}+\frac{9}{8}\right) = 1. \tag{183}$$

This proves part (e). It is interesting to note that had we applied the analogue of Theorem 3 in the $N$-aspect, i.e. directly considered $H_k^*(N)$ without breaking parity, we would only obtain (183) with the inequality $\leqslant 1$ rather than $< 1$. The reason is that as noted in (173), $\alpha_2(O) = 1$. Thus apparently the strict inequality in (183) requires an analysis beyond the diagonal (through breaking parity).

Finally the case (f) of Corollary 2. Again we write

$$\widehat{W(Sp)}(\xi) = \delta_0 + m(\xi).$$

Now $m(\xi)$ is not constant on $(-\frac{4}{3}, \frac{4}{3})$ so that the test function

$$\phi(x) = \left(\frac{\sin 2\pi x}{2\pi x}\right)^2$$

is not optimal. Thus, if

$$\alpha_{\frac{4}{3}}(Sp) := \inf \int_{-\infty}^{\infty} W(Sp)(x)\phi(x)dx, \qquad (184)$$

the infimum being over all $\phi \geqslant 0$, $\phi(0) = 1$ and $\text{support}(\hat{\phi}) \subset (-\frac{4}{3}, \frac{4}{3})$, then

$$\alpha_{\frac{4}{3}}(Sp) < \int_{-\infty}^{\infty} \left(\frac{\sin \frac{4\pi x}{3}}{\frac{4\pi x}{3}}\right)^2 \left(1 - \frac{\sin 2\pi x}{2\pi x}\right) dx = \frac{9}{32}. \qquad (185)$$

Hence

$$\overline{\lim_{K \to \infty}} \frac{1}{M(K)} \sum_{\substack{k \leqslant K \\ k \equiv 0 (\text{mod} 2)}} \sum_{\substack{f \in H_k^*(\Gamma) \\ \text{ord}(\frac{1}{2}, \text{sym}^2(f)) \geqslant 2}} 2 < \frac{9}{32}. \qquad (186)$$

This implies that (recall $L(s, \text{sym}^2(f))$ all have even functional equations)

$$\underline{\lim}_{K \to \infty} \frac{1}{M(K)} \sum_{\substack{k \leqslant K \\ k \equiv 0 (\text{mod} 2)}} \sum_{\substack{f \in H_k^*(\Gamma) \\ \text{ord}(\frac{1}{2}, \text{sym}^2(f)) = 0}} 1 > \frac{55}{64}. \qquad (187)$$

This completes the proof of Corollary 2.

# 5 Appendix. A related extremal problem

In the previous Section we faced particular instances of the following extremization problem. We are given a weight $W(x)$ on $\mathbf{R}$ whose Fourier transform $\widehat{W}(\xi)$ is known only partially, say in the interval $[-2,\ 2]$. The problem is to determine

$$\inf_{\phi} \frac{\int_{-\infty}^{\infty} \phi(x)W(x)dx}{\phi(0)}, \tag{188}$$

the infimum being taken over all $\phi \geq 0$ for which $\text{support}(\hat{\phi}) \subset [-2,\ 2]$. We assume further that $\phi \in L^1(\mathbf{R})$. Examples of $W(x)$ are the densities in (4).

In Section 4 we used the test function

$$\phi(x) = \left(\frac{\sin(2\pi x)}{2\pi x}\right)^2$$

in (188). We show below that it yields almost but not optimal results.

As it stands (188) looks like a linear program problem. However as is pointed out by Gallagher [Ga] it follows from a theorem of Ahiezer and the Paley-Wiener theorem that the admissible functions $\phi$ in (188) coincide with (or have the form)

$$\phi(z) = |h(z)|^2, \tag{189}$$

where $h$ is an entire function of exponential type 1 and $h \in L^2(\mathbf{R})$. That is

$$\hat{\phi}(\xi) = (g * \check{g})(\xi) \tag{190}$$

where

$$\check{g}(\xi) = \overline{g(-\xi)}, \quad \text{support}(g) \subset [-1,\ 1], \quad g \in L^2[-1,\ 1]. \tag{191}$$

By Plancherel Theorem (188) is equivalent to minimization problem

$$\inf_{g \in L^2[-2,\ 2]} R(g),$$

where

$$R(g) = \frac{\int_{-2}^{2} \hat{W}(\xi)(g * \check{g})(\xi)d\xi}{\int_{-2}^{2} (g * \check{g})(\xi)d\xi}. \tag{192}$$

For what we have in mind $\hat{W}(\xi)$ takes the form

$$\hat{W}(\xi) = \delta_0 + m(\xi) \tag{193}$$

for $|\xi| \leq 2$, $m(\xi)$ being a real piecewise continuous function on $[-2, 2]$ which moreover is even in $\xi$. This is the form of the problem which we examine.

Define the self-adjoint operator $K$ from $L^2[-1, 1] \to L^2[-1, 1]$ by

$$Kg(x) = \int_{-1}^{1} m(x - y)g(y)dy. \tag{194}$$

The functional $R$ takes the form

$$R(g) = \frac{<(I + K)g, \, g>}{|<g, \, 1>|^2}. \tag{195}$$

So the minimization problem is that of a quadratic form subject to a linear constraint. Since $W \geq 0$, $R(g) \geq 0$ for any $g$ and so $I + K \geq 0$. Now $K$ is compact and self-adjoint so it has eigenvalues $\lambda_j$, $j = 1, \cdots$, with $|\lambda_j| \to 0$. From the above we have $-1 \leqslant \lambda_j$, for $j = 1, 2, \cdots$.

It may happen that $-1$ is an eigenvalue or equivalently that the finite dimensional kernel

$$N = \ker(I + K) \neq \{0\}.$$

In this case if there is $g \in N$ such that $<g, \, 1> \neq 0$, then clearly $R(g) = 0$, and the minimization in question yields the value 0.

If we are not in this singular case (which will happen if $\lambda_1 > -1$ which is what occurs in our applications), then $\ker(I + K)$ is orthogonal to 1. Hence by Fredholm theory, $1 \in \text{Image}(I + K)$. That is there is an $f_0 \in N^\perp$ such that

$$(I + K)f_0 = 1. \tag{196}$$

Moreover since $I + K > 0$ on $N^\perp$ we have

$$A =< (I + K)f_0, \ f_0 >=< 1, \ f_0 > \tag{197}$$

is positive.

**Proposition 1**: In the nonsingular case

$$\inf_{g \in L^2[-1,\ 1]} R(g) = \frac{1}{A}, \tag{198}$$

and it is attained by $f_0$ which satisfies (196). Moreover if $\lambda_1 > -1$ then this minimizer is unique (in fact the solution to (196) is unique).

Proof: Let $g \in L^2[-1,\ 1]$ with the normalization $< 1,\ g >= A$ (which we can assume if $< 1,\ g > \neq 0$). Then writing $g = f_0 + h$, we have $< h,\ 1 >= 0$, i.e.

$$< h,\ (I + K)f_0 >= 0. \tag{199}$$

Hence

$$\begin{aligned}
R(g) &= \frac{< f_0 + h,\ (I + K)(f_0 + h) >}{A^2} \\
&= \frac{1}{A} + \frac{< h,\ (I + K)(h) >}{A^2} + 2\frac{< h,\ (I + K)(f_0) >}{A^2} \\
&\geqslant \frac{1}{A}.
\end{aligned}$$

So as long as we are not in the singular situation, the minimizer is given by (196), which is a standard Fredholm equation of the second kind. It can be solved in a number of ways.

We can now answer our main question:

**Corollary 2**:

$$\phi(x) = \left(\frac{\sin(2\pi x)}{2\pi x}\right)^2$$

is optimal iff

$$\int_{-1}^{1} m(x-y)dy$$

is independent of $x$.

Proof: If $\phi$ is as the above, the corresponding $g$ is constant on $[-1, 1]$. Now $g = f_0$ is constant according to (196), iff the constant function is an eigenfunction of $K$. This is equivalent to the statement of the corollary.

We apply Corollary 2 to our weights in (4). Firstly, since the property

$$\int_{-\infty}^{\infty} \phi(x)W(x)dx = 0$$

with $\phi \geqslant 0$ and $\phi \in L^1(\mathbf{R})$ implies that $\phi \equiv 0$, it follows that $\lambda_1 > -1$ (in the corresponding eigenvalue problem). That is $(I + K)$ is invertiable. Hence we don't have to worry about the singular case and the unique minimizer $f_0$ satisfies the equation

$$(I+K)f_0 = 1. \tag{200}$$

The functions $m$ for our densities are as follows:

$$m(SO(\text{even}))(\xi) = \frac{1}{2}I_{[-1,\ 1]}(\xi)$$

$$m(SO(\text{odd}))(\xi) = 1 - \frac{1}{2}I_{[-1,\ 1]}(\xi)$$

$$m(Sp)(\xi) = -\frac{1}{2}I_{[-1,\ 1]}(\xi)$$

$$m(O)(\xi) = \frac{1}{2}. \tag{201}$$

According to Corollary 2 in all cases except the last (i.e. $O$), the function $m(\xi)$ is not constant on $[-2, 2]$ and it follows that except in the last case

$$\phi(x) = \left(\frac{\sin(2\pi x)}{2\pi x}\right)^2$$

is not the minimizer. This establishes (171), (172), (173) and (185) of Section 4.

It is not difficult to determine the extremal functions for the $m$'s in (201). J.Vanderkam [Va2] first pointed out to us these functions which he obtained by a direct analysis of the functional (188).

We must solve the equation (196), that is

$$f_0(x) + \int_{-1}^{1} m(x-y) f_0(y) dy = 1, \qquad (202)$$

where $m$ is any one of the functions in (201). Since $m$ is even and the solution $f_0$ is unique, it follows that $f_0$ is an even function of $x$. In particular for $SO(\text{even})$ it satisfies

$$\frac{1}{2} \int_0^1 f_0(y) dy + \frac{1}{2} \int_0^{1-x} f_0(y) dy + f_0(x) = 1, \qquad (203)$$

for $0 \leqslant x \leqslant 1$.

For $SO(\text{odd})$ it satisfies

$$\frac{3}{2} \int_0^1 f_0(y) dy - \frac{1}{2} \int_0^{1-x} f_0(y) dy + f_0(x) = 1, \qquad (204)$$

for $0 \leqslant x \leqslant 1$.

Solving for these by trigonometric functions (i.e. Fourier series) one finds that

$$f_0(x) = \frac{\cos(\frac{x}{2} - \frac{\pi+1}{4})}{\sqrt{2}\sin(\frac{1}{4}) + \sin(\frac{\pi+1}{4})}, \quad 0 \leqslant x \leqslant 1, \qquad (205)$$

solves (203), while its even extension to $[-1, 1]$ solves (202). Applying (197) and Proposition 1 then yields the minimum. A calculation then shows that

$$\alpha_2(SO(\text{even})) = \frac{3 + \cot(\frac{1}{4})}{8} = 0.8645\cdots. \qquad (206)$$

Similarly

$$f_0(x) = \frac{\cos(\frac{x}{2} + \frac{\pi-1}{4})}{3\sin(\frac{\pi+1}{4}) - 2\sin(\frac{\pi-1}{4})}, \qquad (207)$$

for $0 \leqslant x \leqslant 1$, solves (204). It leads to the value

$$\alpha_2(SO(\text{odd})) = \frac{5 + \cot(\frac{1}{4})}{8} = 1.1145\cdots. \tag{208}$$

In particular these allow us to conclude that

(a). (176) of Section 4 holds with $> \frac{9}{16}$ replaced by

$$\geqslant \frac{13 - \cot(\frac{1}{4})}{16} = 0.5678\cdots$$

(b). (180) of Section 4 holds with $> \frac{15}{16}$ replaced by

$$\geqslant \frac{19 - \cot(\frac{1}{4})}{16} = 0.94275\cdots$$

(c). (183) of Section 4 holds with $< 1$ replaced by

$$\leqslant \frac{4 + \cot(\frac{1}{4})}{8} = 0.9895\cdots$$

**Acknowledgement**: We would like to thank S.J.Miller for his careful reading and comments on this paper.